\documentclass[twoside]{irmaems}
\usepackage{amssymb} %use \usepackage[psamsfonts]{amssymb} if your fonts are from Y&Y/Blue Sky Research
\usepackage{amsmath} %use \usepackage{amsmath,cmmib57} if your fonts are from Y&Y/Blue Sky Research
\usepackage{latexsym}

\renewcommand\theequation{\arabic{section}.\arabic{equation}}

\theoremstyle{definition} %%% for statements in roman typeface

 \newtheorem{definition}{Definition}[section]

\theoremstyle{plain}      %%% for statements in italic typeface

 \newtheorem{theorem}[definition]{Theorem}
 \newtheorem{corollary}[definition]{Corollary}
 \newtheorem{lemma}[definition]{Lemma}

\newtheorem*{conjecture}{Conjecture}

\begin{document}

\title{The Theory of Quasiconformal Mappings in Higher Dimensions, {\bf I}}

\author{Gaven J. Martin\thanks{
Work partially supported by the New Zealand Marsden Fund}
\address{
Institute for Advanced Study\\
Massey University,  Auckland\\
New Zealand\\
email:\,\tt{g.j.martin@massey.ac.nz}
\\
}}

\maketitle

\begin{abstract} This chapter presents a survey of the many and various elements of the modern higher-dimensional theory of quasiconformal mappings and their wide and varied application.  It is unified (and limited) by the theme of the author's interests.  Thus we will discuss the basic theory as it developed in the 1960s in the early work of F.W. Gehring and Yu G. Reshetnyak and subsequently explore the   connections with geometric function theory,  nonlinear partial differential equations,  differential and geometric topology and dynamics as they ensued over the following decades.  We give few proofs as we try to outline the major results of the area and current research themes.  We do not strive to present these results in maximal generality,  as to achieve this considerable technical knowledge would be necessary of the reader.  We have tried to give a feel of where the area is,  what are the central ideas and problems and where are the major current interactions with researchers in other areas.  We have also added a bit of history here and there.  We have not been able to cover the many recent advances generalising the theory to mappings of finite distortion and to degenerate elliptic Beltrami systems which connects the theory closely with the calculus of variations and nonlinear elasticity,  nonlinear Hodge theory and related areas,  although the reader may see shadows of this aspect in parts.  

In the sequel (with Bruce Palka) we will give a more detailed account of the basic techniques and how they are used with a view to providing tools for researchers who may come in contact with higher-dimensional quasiconformal mappings from their own area.
\end{abstract}


\newpage

\tableofcontents   

\section{Introduction}  
  
Geometric Function Theory in higher dimensions is largely
concerned with generalisations to ${\Bbb R}^n$,  $n\geq 3$,  of aspects of complex analysis,  the theory of
analytic functions and conformal mappings - particularly the geometric and function theoretic properties.  In this sense it has been a successful theory with a variety of applications,  many of which we will discuss below.  

\medskip

The
category of maps that one usually considers in the higher-dimensional theory are the {\em quasiregular} mappings, 
or,  if injective,  {\em quasiconformal}  mappings.  Both kinds of mappings have the
characteristic property of {\em bounded distortion}.    The higher-dimensional theory of quasiconformal mappings was initiated in earnest by Yu G. Reshetnyak (USSR),  F.W. Gehring (USA) and J. V\"ais\"al\"a (Finland) in 1960-62,  see \cite{Gehring1,Gehring2,Reshetnyak1,Reshetnyak2}.  There was earlier work,  notably that of  Ahlfors-Beurling (1950) on conformal invariants, Ahlfors (1954) and Callender (1959).  While Ahlfors' work was focussed on two-dimensions,  the geometric ideas and techniques had clear generalities.  Callender followed Finn \& Serrin to establish H\"older continuity estimates for higher-dimensional quasiconformal mappings.  We note that one of the most famous applications of the theory of higher-dimensional quasiconformal mappings,  Mostow's rigidity theorem (1967) \cite{Mostow}, came just five years after the basic foundations were laid.   

The generalisations to non-injective mappings was initiated with Reshetnyak and the basic theory was comprehensively laid and significantly advanced in a sequence of papers from the Finnish school of  O. Martio, S. Rickman and J. V\"ais\"al\"a in the late 1960s \cite{MRV1,MRV2,MRV3}.

 \medskip

Both quasiconformal and quasiregular mappings solve natural partial differential equations (PDE) closely analogous to the familiar Cauchy--Riemann and Beltrami
equations of the plane.  The primary difference being that in higher dimensions these equations necessarily become nonlinear and overdetermined.   Other desirable properties for a theory of the geometry of mappings are that they should preserve the natural Sobolev spaces which arise in consideration of the function
theory and PDEs on subdomains of  ${\Bbb R}^n$,  or more generally $n$--manifolds.  Quasiconformal mappings do have these properties.

\medskip

In dimensions $n\geq 3$,  one needs to move away from the class of conformal mappings because of their remarkable rigidity properties.  Perhaps most well known is the Liouville theorem,  basically established in the 1970s independently and by different methods by Gehring and Reshetnyak as we discuss later.  This rigidity is partly explained by the fact that the governing equations are overdetermined.  This rigidity also has consequences for the pronounced differences between injective and non-injective mappings in higher dimensions.  In two-dimensions,  the celebrated Stoilow factorisation theorem asserts that a quasiregular mapping $f:\Omega\to{\Bbb C}$ admits a factorisation $f = \varphi\circ g$ where $g:\Omega\to{\Bbb C}$ is quasiconformal and $\varphi$ is holomorphic.  This factorisation, together with our more or less complete understanding of the structure of holomorphic functions in the plane, connects quasiregular \& quasiconformal mappings strongly.  In particular, if $B_f$ is the branch set,  
\begin{equation}
B_f = \{ x \in \Omega : \mbox{ $f$ is not locally injective at $x$} \}, 
\end{equation}
then, in two-dimensions, the factorisation theorem quickly shows that $B_f$ will be a discrete subset of $\Omega$.  Thus if $\Omega$ is the Riemann sphere, ${\hat{{\Bbb C}}}={\Bbb C}\cup\{\infty\}\approx {\Bbb S}^2$, $B_f$ will be a finite point set.  This is far from true in higher dimensions.  First,  well known results in geometric measure theory,  see Federer \cite{Fed},  connect smoothness (in terms of differentiability) and local injectivity. Thus branched maps cannot be very smooth ($C^1$ is alright in three dimensions though).  Second,  the branch set $B_f$ of a mapping of bounded distortion can have quite pathological topology,  for instance it could be Antoine's necklace - a Cantor set in ${\Bbb S}^3$ whose complement is not simply connected.  This makes Rickman's development of the higher-dimensional Nevanlinna theory all the more remarkable.   As a consequence of the Nevanlinna theory,  as in the classical case,  one obtains best possible results concerning precompactness of families of mappings.  Since quasiregular mappings are open and discrete at this point it is work recalling Chernavskii's theorem \cite{Chern,Vaisala1} which asserts that if $B_f$ is the branch set of a quasiregular mapping,  then  the topological dimension of both $B_f$ and $f(B_f)$ is less than or equal to $n-2$ and therefore cannot separate.  Further, fairly general topological results enable one to talk about the degree and topological index of such mappings.
\medskip

There are also  second order equations related to the nonlinear governing equations for quasiconformal mappings.  For example,  the components
of an analytic function are harmonic,  while those of a quasiregular mapping are  ${\cal A}$--harmonic.  These are basically the Euler-Lagrange equations for a conformally invariant integral for which the mapping in question is a minimiser.  In this way such well--known non-linear differential operators as the $p$--Laplacian and
the associated non-linear potential theory arise naturally in the theory of higher-dimensional quasiconformal mappings.  This potential theory has significant topological implications for mappings of bounded distortion.  These were first observed by Reshetnyak. 

\medskip

Another fruitful idea when studying quasiconformal mappings and their properties is to view quasiregular
mappings as conformal with respect to certain measurable Riemannian or measurable conformal structures.   In two-dimensions this gives the direct connection with Teichm\"uller theory of course and this idea is greatly aided by the fact one can solve the associated Beltrami equation,  leading to the so called measurable Riemann mapping theorem - or the existence theorem for quasiconformal mappings.   In higher dimensions,  unfortunately almost nothing useful is known about solving Beltrami systems.  There are obvious reworkings of the classical results from the 1920s of Weyl and Schouten,  which assume the vanishing of a second order tensor,  when the conformal tensor  $G$ is sufficiently smooth. It is an extremely interesting problem to try and give reasonable conditions on $G$ which guarantee local existence if $G$ is perhaps only $C^{1+\epsilon}$ smooth let alone the most important case when $G$ is only assumed measurable. The general higher-dimensional theory does provide good results about the regularity of solutions,  really initiated by Gehring's  higher integrability and reverse H\"older inequalities from 1973 \cite{GehHI}.  These generalise earlier results of Bojarski from 1955 in two-dimensions \cite{B}, but again totally different methods are needed to attack these non-linear equations in higher dimensions.  I think it is fair to say the Gehring's higher integrability results revolutionised the theory and ultimately brought it closer to PDE and nonlinear analysis as the techniques he developed had much wider application.  We also understand,  to a reasonable extent at least,  both the uniqueness of solutions to higher-dimensional Beltrami systems as well as analytic continuation and so forth.  

\medskip

Quasiconformal mappings provide  a class of deformations which lie ``between'' homeomorphisms and  diffeomorphisms but enjoy compactness properties neither do.  The most recent developments in the theory concern mappings of {\em finite distortion}.  Here the assumption concerning boundedness is removed and replaced with various control assumptions on the distortion or its associated tensors.  Such mappings are even more flexible and to study them more refined techniques are necessary as the governing equations will be degenerate elliptic.  However various compactness properties of families of
mappings with finite distortion make them ideal tools for solving various problems in
$n$-dimensional analysis.  For instance in studying deformations of
elastic bodies and the related extremals for variational integrals,  mappings of finite
distortion are often the natural candidates to consider.  These ideas lead directly to the theory of non-linear elasticity  developed by Antman and Ball,  and many others.  This theory of elasticity studies mappings (in certain
Sobolev classes) which minimize various stored energy integrals.  On seeks existence,  regularity and so forth.  The Jacobian determinant,  in particular,  has been subjected to a
great deal of investigation.   

\medskip

Of course there are many outstanding problems which are helping to  drive the field,  but which we won't discuss here.  These include determining precise geometric conditions on a domain to be quasiconformally equivalent to a ball (thus a generalised Riemann mapping theorem).  As we will see in a moment,  the Liouville theorem implies that in dimensions $n\geq 3$ any domain conformally equivalent to the unit ball is a round ball or half-space.  In two-dimensions Ahlfors gave a beautiful intrinsic characterisation of the quasiconformal images of the unit disk.  While such a nice result is unlikely in higher dimensions, not a great deal is really known.  

In other directions,  Iwaniec and his coauthors are advancing the connections between the higher-dimensional theory (largely as it pertains to the geometry of mappings)  and the calculus of variations.  Particular advances concern generalising the theory to mappings of finite distortion \cite{IKM,AIKM}.  Here the distortion is no longer assumed uniformly bounded,  but some additional regularity is necessary to get a viable theory.  Most of the major results assume something close to the distortion function $K(x,f)$ being of bounded mean oscillation.   As the Jacobian of a mapping $J(x,f)$ has automatically higher regularity as an $H^1$ function,  one seeks to exploit the $H^1$-BMO duality discovered by Fefferman to gain information about the total differential from the distortion inequality
\[ |Df(x)|^n \leq K(x,f) J(x,f) \]
since meaning can be given to the right-hand side.  There are very many interesting problems and deep connections to other areas here.

There are still further generalisations and applications in the geometry and analysis of metric spaces.  The connections with the higher-dimensional theory of quasiconformal mappings was pioneered by Heinonen and Koskela \cite{HK} and is a very active area of research today.

\medskip

Here is a sample of the successful and diverse applications of the higher-dimensional theory of quasiconformal mappings (some mentioned above and in no particular order of importance): 

\begin{itemize}
\item Compactness,  equicontinuity and local to global distortion estimates;  
\item The Liouville theorem and other stability and rigidity phenomena; 
\item Gehring's improved regularity and higher integrabilty; 
\item Mostow rigidity - uniqueness of hyperbolic structures ($n\geq 3$); 
\item Sullivan's uniformisation theorem - the existence of quasiconformal structures on topological $n$-manifolds ($n\neq 4$); 
\item Rickman's versions of the Picard theorem and Nevanlinna theory; 
\item Applications in nonlinear potential theory, $A$-harmonic functions and non-linear elasticity; 
\item Tukia-V\"ais\"al\"a's ``quasiconformal geometric topology''; 
\item Quasiconformal group actions and geometric group theory; 
\item Donaldson and Sullivan's ``quasiconformal Yang-Mills theory''; 
\item Painlev\'e type theorems and the structure of singularities; 
\item Quasiconformal maps in metric spaces with controlled geometry;
\item Analysis and geometric measure theory in metric spaces.
\end{itemize}

\medskip

Mindful of the readership of a chapter such as this,  we will not strive for maximum generality in the results we present.  Also,  we will seldom present complete proofs and discussions,  but set the reader toward places where such discussions can be found.  To that end there are a number of relatively recent books which the reader might consult for details of omissions here and which have a relatively broad focus.  A reasonably complete account of the modern two-dimensional theory is given in  Astala, Iwaniec \& Martin, {\em Elliptic partial differential equations and quasiconformal mappings in the plane}, 2009,  \cite{AIM}.   For the analytic aspects of the theory in higher dimensions we have Iwaniec \& Martin, {\em Geometric function theory and non-linear analysis}, 2001, \cite{IM}.  For the nonlinear potential theory see Heinonen,  Kilpel\"ainen, \& Martio, {\em Nonlinear potential theory of degenerate elliptic equations}, Oxford, 1993, \cite{HKM}.  For analysis on metric spaces see Heinonen, {\em Lectures on analysis on metric spaces},  2001 \cite{H}.  For the Nevanlinna and related theories of quasiregular mappings we have Rickman, {\em Quasiregular mappings},   1993 \cite{R}.  Vuorinen, {\em Conformal geometry and quasiregular mappings},  1988, \cite{Vu} gives a detailed account of the distortion estimates and other geometric aspects of the theory which is further developed in by  Anderson,   Vamanamurthy and  Vuorinen,  \cite{AVV}.

Of course there are others,  but these books should give a more or less complete overview.  But not to forget the past,  we cannot fail to mention the classics,  Ahlfors, {\em Lectures on quasiconformal mappings},  1966, \cite{A}  and Lehto \& Virtanen, {\em Quasiconformal mappings in the plane}, 1973, \cite{LV}, for the two-dimensional theory and of course V\"ais\"al\"a {\em Lectures on $n$-dimensional quasiconformal mappings},  1971, \cite{Vaisala2},  a book from which many of us learnt the basics of the higher-dimensional theory.

\section{Two geometric definitions.}

We will present the analytic definition of a quasiconformal mapping via Sobolev spaces and differential inequalities a little later.  However we want to give a brief initial discussion to capture the idea of infinitesimal distortion.  This is because the geometric definitions of quasiconformality are quite global in nature - asking us to test a Lipschitz condition against every family of curves in a given domain.  It is this interplay between the local definition of quasiconformality and the global one that is a real strength of the theory.  Once one has established an infinitesimal distortion condition (through properties of solutions to a PDE or some assumptions around differentiability), then one obtains large scale distortion estimates through considering various curve families and geometric estimates upon them.

\subsection{The linear distortion}

Let $\Omega$ and $\Omega^\prime$ be domains in ${\Bbb R}^n$ and let $f:\Omega \to \Omega^\prime$ be a
homeomorphism.  We will define quasiconformal\index{quasiconformal!linear dilatation definition}
mappings as mappings of ``bounded distortion'' and therefore we must discuss what distortion might mean.  Suppose therefore, that $x\in \Omega$ and $r<d(x,\partial \Omega)$.  We define the infinitesimal distortion $H(x,f)$ of $f$ at the point $x$ as 
\begin{equation}\label{oH} 
H(x,f) = \limsup_{|h|\to 0}\;\; \frac{\max |f(x+h)-f(x)|}{\min |f(x+h)-f(x)|}
\end{equation} 
We further say that $f$ is {\em quasiconformal} in $\Omega$ if
$H(x,f)$ is bounded throughout the domain 
$\Omega$:
\[ \mbox{\em there exists an $H< \infty$ such that $H(x,f) \leq H$ for every $x\in \Omega$} \] 
The essential supremum of this quantity is called the {\em linear distortion}\index{linear distortion} of $f$.   
\begin{equation}
H(f) = \| H(x,f)\|_{L^\infty(\Omega)}.
\end{equation}
Notice the pointwise {\em everywhere} assumption here in the definition of quasiconformality.  It is necessary.
If $f$ is differentiable at $x_0$ with non-singular differential, then we can multiply and divide by $|h|$ in (\ref{oH}) and   take limits.  Quasiconformality  quickly yields an inequality between the smallest and largest directional derivatives, 
\[ \max_{|h|=1} |Df(x_0) h| \leq H(x_0,f) \; \min_{|h|=1} |Df(x_0) h|. \] 
A little linear algebra reveals that the left hand side here is the square root of the largest eigenvalue of the positive definite matrix $D^tf(x_0)Df(x_0)$,  and the right hand side is the smallest such.

While geometrically appealing,  unfortunately this
quantity is not particularly useful in higher dimensions since it is not lower semicontinous on the space of
quasiconformal mappings mappings \cite{Ifail}.  Recall that a real-valued function $h$ is {\em lower semicontinuous}
\index{lower semicontinuous} if for all $x_0$, 
\[ \liminf_{x\to x_0} h(x) \geq h(x_0). \]
A lower semicontinuous distortion function  will guarantee distortion does not suddenly increase in the limit, a clearly
desirable property.  Interestingly,  this failure is directly connected with the failure of rank-one convexity in the calculus of variations.  However,  there is a remarkable result here due to Heinonen and Koskela \cite{HK2}.  It turns out the the $\limsup$ requirement of the definition at (\ref{oH}) is met once a $\liminf$ condition holds.\index{quasiconformal!$\liminf$ definition}

\begin{theorem}  Let $f:\Omega\to{\Bbb R}^n$ be a homeomorphism.  If there is $H^*<\infty$ such that for every $x\in \Omega$ 
\begin{equation}\label{liH} 
H^*(x,f) = \liminf_{|h|\to 0}\;\; \frac{\max |f(x+h)-f(x)|}{\min |f(x+h)-f(x)|} \leq  H^*,
\end{equation}
then there is $H=H(n,H^*)<\infty$ such that for every $x\in \Omega$
\begin{equation} 
H(x,f) \leq  H,
\end{equation}
and consequently (\ref{liH}) is enough to guarantee the quasiconformality of $f$.
\end{theorem}
Again,  note the requirement of having a condition at every point of $\Omega$.  The analytic definitions of quasiconformality will get around this problem by having a pointwise almost everywhere criteria.  However, these conditions  must of course give the boundedness of the linear distortion everywhere. Before we go in that direction we discuss the earliest natural definition of quasiconformality which is through the bounded distortion of a  conformal invariant -- the moduli of curve families.

\subsection{Moduli of curve families}

 The most useful geometric definition  of a quasiconformal mapping is through the modulus of a curve family.\index{curve family!definition}  A curve family
${\Gamma}$ is simply a collection of (rectifiable) curves;  continuous maps $\gamma:[0,1]\to\Omega$.  It is usual to identify a curve with its image as the quantities we wish to study must be invariant of parameterisation.  Typically a curve family will be of the following sort, $\Delta(E,F:\Omega)$ the set of all curves connecting $E$ to $F$ and lying in a domain  $\Omega$ of ${\Bbb R}^n$. 
Given a curve family
${\Gamma}$, an {\em admissible density}\index{admissible density} is a non-negative Borel function $\rho$ for which
\begin{equation}\label{modulusdefn}
\int_\gamma \rho(s) \; ds \geq 1, \hskip20pt  \mbox{ for all } \gamma \in {\Gamma}.
\end{equation}
We see immediately that highly irregular curves (in particular those that are not rectifiable) in a family will not be relevant as their $\rho$-length will most likely be infinite,  while constant curves cannot admit such a $\rho$.

\medskip

The {\em modulus} \index{curve family!modulus} of ${\Gamma}$ is
\begin{equation}
M({\Gamma}) = \inf \int_{{\Bbb R}^n} \rho^n(x)\; dx
\end{equation}
where the infimum is over all admissible densities for ${\Gamma}$.  

\medskip

There are more general discussions to be had here.  We could,  for instance, consider the $p$-modulus (or $p$-capacity) where we look at 
\[ \inf \int_{{\Bbb R}^n} \rho^p(x)\; dx\]
over the same class of admissible functions (this was first considered by Fuglede).  These quantities can be used to detect the size of sets in a similar fashion to Hausdorff dimension.  As an example consider a set $E\subset {\Bbb B}^n$ and let ${\Gamma}$ consist of all the curves in ${\Bbb B}^n$ connecting $E$ to ${\Bbb S}^n$.  Then, depending on the size and structure of $E$, there may be a value $p_0$ for which this quantity is zero -  $E$ has $p_0$-capacity zero.  This has geometric and function theoretic consequences.  For instance,  sets of $n$-capacity (usually called conformal capacity) zero are typically negligible in the theory of quasiconformal mappings and so, for instance,   removable for bounded mappings and so forth.  However,  these sets have Hausdorff dimension zero and so are very thin.

\bigskip

The idea of the modulus of curve families is to develop the ``length-area'' method used by Ahlfors and Beurling to great effect in two-dimensions in their celebrated paper on conformal invariants  \cite{AB} in 1950,  although these ideas had been around and used in various ways in complex analysis since the 1920's. 

\bigskip

There are a few basic properties of the modulus of a curve
family which fall out of the definition.  Firstly $M({\Gamma})$ is increasing. If ${\Gamma}_1\subset{\Gamma}_2$,  then $M({\Gamma}_1)
\leq M({\Gamma}_2)$.  If ${\Gamma}$ contains a single ``constant curve'',  then $M({\Gamma})=\infty$.  If ${\Gamma}_1$ and ${\Gamma}_2$
are curve families such that every curve in ${\Gamma}_2$ has a subcurve in ${\Gamma}_1$,  then $M({\Gamma}_2)\leq M({\Gamma}_1)$.
Finally  modulus is countably subadditive and additive on disjoint families. 

The most important fact is of course that modulus is a conformal invariant. For the moment a conformal
mapping  will be a diffeomorphism whose differential (the matrix $Df$) is pointwise a scalar multiple of an
orthogonal matrix. For a conformal map we therefore have the equality
\begin{equation}\label{confdef}
|Df(x)|^n=J(x,f) = \det Df(x)
\end{equation}
(recall $|Df(x)|=\max_{|h|=1}|Df(x)h|$,  the largest directional derivative).

\medskip

Next we establish that the modulus is conformally invariant.  The proof is easy,  but the reader should take note of how the differential inequalities between $|Df(x)|$ and $J(x,f)$ are used as this motivates the analytic definition. 

\begin{theorem}\label{mci}
Modulus is a conformal invariant.
\end{theorem}
\noindent{\bf Proof.} Let ${\Gamma}$ be a curve family and set ${\Gamma}^\prime=f({\Gamma})$.  If $\rho_1$ is
admissible for
${\Gamma}$,  then
$\rho(x)=\rho_1(f(x))|Df(x)|$ is admissible for ${\Gamma}$ since
\begin{equation}
\int_\gamma \rho \; ds = \int_\gamma \rho_1(f(x))|Df(x)|\; ds = \int_{\gamma^\prime} \rho_1 \; ds \geq 1
\end{equation}
Next  
\begin{eqnarray*}
M({\Gamma}) & \leq & \int_{{\Bbb R}^n} \rho^n(x)\; dx = \int_{{\Bbb R}^n} \rho_{1}^{n}(f(x)) |Df(x)|^n\;  dx \\
& = & \int_{{\Bbb R}^n} \rho_{1}^{n}(f(x)) J(x,f)\; dx = \int_{{\Bbb R}^n} \rho_{1}^{n}(x)\; dx .
\end{eqnarray*}
Now taking the infimum over all $\rho_1$ shows us that  $M({\Gamma}) \leq M({\Gamma}^\prime)$.  The converse inequality holds
since $f^{-1}$ is also a conformal mapping. \hfill $\Box$

\medskip

We can now give an alternative definition of quasiconformality\index{quasiconformal!geometric definition}

\subsection{The geometric definition of quasiconformality} 
Let $f:\Omega \to \Omega^\prime$ be a homeomorphism.  Then $f$ is $K$--quasiconformal if there exists a $K$, $1 \leq K \leq
\infty$,  such that
\begin{equation}\label{Moddist}
\frac{1}{K} M({\Gamma}) \leq M(f{\Gamma}) \leq K M({\Gamma})
\end{equation}
for every curve family ${\Gamma}$ in $\Omega$.

\medskip

Of course in practise it is impossible to test the condition (\ref{Moddist}) against {\em every} curve family.  That is why we seek equivalent conditions - either infinitesimal or by testing against only certain curve families - which guarantee quasiconformality.  Once we have such things at hand (\ref{Moddist})  provides powerful global geometric information - provided we can find ways of computing,  or at least estimating,  the moduli of curve families.

\medskip

There are a couple of direct consequences from this definition that are not nearly so trivial when using the analytic definitions that follow.  Primary among these are

\begin{theorem}\label{composition}  Let $f:\Omega\to\Omega'=f(\Omega)$ be $K$-quasiconformal and $g:\Omega'\to{\Bbb R}^n$ be $K'$-quasiconformal.  Then
\begin{itemize}
\item  $f^{-1}:\Omega'\to\Omega$ is $K$-quasiconformal,  and 
\item $g\circ f:\Omega \to {\Bbb R}^n$ is $KK'$-quasiconformal.
\end{itemize}
\end{theorem}

The general theory now develops by computing the modulus of special sorts of curve families.  Then we estimate the modulus of more general curve families in terms of geometric data and define various special functions for the
modulus of various curve families which are in some sense extremal for moduli problems (for instance the Gr\"otzsch 
and Teichm\"uller curve families being the most common).  We then obtain geometric information about quasiconformal mappings by studying
what happens to special curve families under quasiconformal mappings using the Lipschitz estimates at  (\ref{Moddist}) and
comparing with the general estimates of various moduli.  This approach quite quickly reveals that quasiconformal
mappings are locally H\"older continuous,  and establishes such things as an appropriate version of the Schwarz Lemma, 
and so forth. 

Of course the Lipschitz estimate at (\ref{Moddist}) must have further consequences for the differentiability and regularity of the homeomorphism $f$.  A major part of the basic theory is in identifying these. This is typically done by connecting this geometric definition,  with the analytic definition we give a bit later.   To get further into the theory we must actually compute a couple of examples of the modulus of curve families.

\subsection{The modulus of some curve families}

First,  and most useful,  is the calculation of the modulus of the curves in an annular ring.
\begin{theorem}
Let ${\Gamma} = \Delta({\Bbb S}^{n-1}(a),{\Bbb S}^{n-1}(b):{\Bbb B}^n(b)\setminus \overline{{\Bbb B}^n(a)})$, the set of curves connecting the boundary components of the annulus $A(a,b)={\Bbb B}^n(b)\setminus \overline{{\Bbb B}^n(a)}$.  Then, 
\begin{equation} 
M({\Gamma}) =    \omega_{n-1} \biggl(\log \; \frac{b}{a} \biggr)^{1-n}
\end{equation}  
where $\omega_{n-1}$ is the volume of the $(n-1)$-dimensional sphere.
\end{theorem}
\noindent{\bf Proof.}  Choose an  admissible density $\rho$ for ${\Gamma}$. The rays  $\gamma_u (r) = ru$, $u\in {\Bbb S}^{n-1}$ and $a<r<b$ lie in ${\Gamma}$ and so H\"older's inequality gives us 
\begin{eqnarray*}
1 & \le & \biggl( \int_{\gamma_u} \rho\, ds\biggr)^n 
= \biggl( \int_a^b \rho (ru)\, dr\biggr)^n = \biggl( \int_a^b \rho (ru)\, r^{(n-1)/n} r^{(1-n)/n} \: dr\biggr)^n \\
& = &\biggl( \int_a^b \left[\rho (ru)\, r^{(n-1)/n} \right]^n \biggr) \biggl(\int_a^b r^{-1} \: dr\biggr)^{n-1}\\ &\leq  & \biggl( \int_a^b \rho^n (ru)\,r^{n-1}\, dr\biggr) 
\biggl( \log \frac{b}{a} \biggr)^{n-1} .
\end{eqnarray*}
Because this  holds for every $u$ in $S^{n-1}$, we can integrate it   
over $S^{n-1}$.  Thus
\begin{eqnarray*}
\int_{{\Bbb R}^n} \rho^n\, dx  \geq \int_A \rho^n \, dx = \int_a^b \biggl[ \int_{{\Bbb S}^{n-1}} \rho^n (ru)\, d\sigma_{n-1} 
\biggr] r^{n-1}\, dr \\
= \int_{{\Bbb S}^{n-1}} \biggl[ \int_a^b \rho^n (ru)\, r^{n-1}\, dr\biggr] \, 
d\sigma  \geq   \omega_{n-1} \biggl(\log \; \frac{b}{a} \biggr)^{1-n}
\end{eqnarray*}
giving us the lower bound we want.  Next, define $\rho :{\Bbb R}^{n} \to [0,\infty]$  by 
\[ \rho (x) = \begin{cases} \big(\log \frac{b}{a}\big)^{-1}\,  |x|^{-1} &\mbox{if $x\in A$\ ,}\\ 
0&\mbox{if $x\notin A$\ .}
\end{cases}
\] 
Then $\rho$ is an admissible function since
\[ \int_\gamma\rho\, ds \geq \int_a^b \rho (r)\, dr = \int_a^b \big(\log \frac{b}{a}\big)^{-1}\,  r^{-1}\, dr= 1\] 
for every  rectifiable $\gamma$ in ${\Gamma}$.  Then
\begin{eqnarray*}
M({\Gamma})  &\leq& \int_{{\Bbb R}^n} \rho^n \, dx
=  \int_a^b \biggl[ \int_{{\Bbb S}^{n-1}} \rho^n (ru)\, d\sigma \biggr] 
r^{n-1}\, dr \\ & = &   \frac{1}{\log(b/a)} \int_a^b \biggl[\int_{{\Bbb S}^{n-1}} r^{-n} \, d\sigma 
\biggr] r^{n-1} \, dr   
= \big(\log\;\frac{b}{a}\big)^{-n} \omega_{n-1} \int_a^b r^{-1}\, dr \\
& = &  \omega_{n-1} \big(\log\;\frac{b}{a}\big)^{1-n}.
\end{eqnarray*}
So we obtain the desired equality.  \hfill $\Box$

\bigskip
 
Unfortunately,  very few other moduli can be explicitly computed.  An elementary estimate is the following

\begin{lemma} If $E$ is an open set in ${\Bbb R}^n$ with whose Lebesgue measure is finite, $|E|<\infty$,  and if ${\Gamma}$ is a family 
of curves in $E$ for which $d= \inf\{ \ell (\gamma) :\gamma\in{\Gamma}\} >0$, 
then 
\[ M  ({\Gamma}) \leq \frac{|E|}{d^n} < \infty\ .\] 
\end{lemma}
This is a direct consequence of the fact that 
$\rho = d^{-1}\chi_E$ is an admissible density  
for ${\Gamma}$.

 \medskip
 
 One now seeks ways to estimate, both from above and below, the modulus of certain curve families.  With experience we quickly find that we are most often concerned with curve families that join two components of the boundary of a domain, and these are called ``rings'',  or sometimes condensers.  Thus $R(E,F;G) $ is the family of curves joining $E$ to $F$ in the domain $G$ 
 and when $G={\Bbb R}^n$ we simply write $R(E,F)$.   The modulus of a ring $R=R(E,F;G)$ (or the capacity of the  condenser) is
 \begin{equation}
 {\rm Mod}(R) = \left[\frac{1}{\omega_{n-1}} M({\Gamma}_R) \right]^{n-1}
 \end{equation}
where ${\Gamma}_R$ is the family of all curves joining $E$ to $F$ in $G$.  In particular
\begin{lemma}  The modulus of the annular ring $A=\{x:a<|x|<b\}$ is
\[ {\rm Mod}(A) = \log\, \frac{b}{a} \]
\end{lemma}
Note now that bigger rings have smaller modulus.  As the components get closer together we can expect the modulus to tend to $\infty$.

Now,  with this notion the Lipschitz estimate of  (\ref{Moddist}) gives rise to a new characterisation of quasiconformality.  Suppose that $f:\Omega\to{\Bbb R}^n$ is a mapping and  suppose that $E,F\subset \Omega$ are continua with $R=R(E,F;\Omega)$.  We define
\[ f(R) = R(f(E),f(F),f(\Omega)).\]
Then we may say that $f$ is quasiconformal if  
\begin{equation}\label{Ring}
\frac{1}{K} {\rm Mod}(R) \leq {\rm Mod}(f(R)) \leq K {\rm Mod}(R)
\end{equation}
for every ring in $\Omega$.  Notice that this requires
\[ \frac{1}{K} \left[\frac{1}{\omega_{n-1}} M({\Gamma}_R) \right]^{n-1}  \leq   \left[\frac{1}{\omega_{n-1}} M({\Gamma}_{f(R)}) \right]^{n-1} \leq    K \left[\frac{1}{\omega_{n-1}} M({\Gamma}_R) \right]^{n-1} \]
and hence
\[ \frac{1}{K^{1/(n-1)} } M({\Gamma}_R)   \leq   M({\Gamma}_{f(R)})  \leq    K^{1/(n-1)} M({\Gamma}_R) \] 
and the precise measure of distortion (namely $K$) differs from that at (\ref{Moddist}) unless $n=2$.  

\medskip

We note this as a warning,  it is not really of relevance unless one is seeking optimal constants regarding various sorts of estimates - continuity and so forth.   Nevertheless,  there is room for confusion.  

\medskip

Aspects of the higher-dimensional theory are concerned with the continuity of the modulus of rings in the Hausdorff topology.  We recall that subcontinua of ${\Bbb R}^n$,  $E_j$,  converge to $E$ in the Hausdorff topology if
\[ \sup_{x\in E} \{{\rm dist}_q(x,E_j)\} +  \sup_{y\in E_j} \{{\rm dist}_q(y,E)\}  \to 0,  \hskip10pt\mbox{as $j\to\infty$} \]
Here ${\rm dist}_q$ refers to the spherical distance of ${\Bbb R}^n$,
\begin{equation}
{\rm dist}_q(x,y) = \frac{|x-y|}{\sqrt{|x|^2+1}\, \sqrt{|y|^2+1}} .
\end{equation}
\begin{theorem}\label{Rcont}  Suppose that $E$ and $F$ are disjoint continua (${\rm dist}_q(E,F)>0$) and that $E_j\to E$ and $F_j\to F$ in the Hausdorff topology.  Then
\[ {\rm Mod}(R(E_j,F_j:{\Bbb R}^n) )\to {\rm Mod}(R(E,F:{\Bbb R}^n)), \hskip20pt \mbox{ as $j\to\infty$} \]
\end{theorem}
Also as the modulus of curve families decreases under inclusion, the modulus of rings increases.
\begin{lemma}  Suppose $E_0\subset E_1$ and $F_0\subset F_1$ in $\Omega$, then
\begin{equation}\label{modincr}
{\rm mod}(R(E_0,F_0:\Omega)) \leq {\rm mod}(R(E_1,F_1:\Omega)).
\end{equation}
\end{lemma}

\subsection{The Gr\"otzsch and Teichm\"uller rings}

We have computed the modulus of the annulus above. What we need now are some more general rings whose modulus we can estimate well and prove some extremal properties for.  
The first is the Gr\"otzsch ring\index{Gr\"otzsch ring}.  We denote for $t>1$
\[ R_G(t) = R({\Bbb B}^n, [t,\infty) :{\Bbb R}^n) \]
where by $ [t,\infty)$ we mean $\{(s,0,\ldots,0): t\leq s \}$.
Next is the Teichm\"uller ring\index{Teichm\"uller ring}.  Here
\[ R_T(t) = R([-1,0], [t,\infty):{\Bbb R}^n) \]
and then we set
\[ \gamma_n(t)= {\rm Mod}(R_G(t)), \hskip20pt \tau_n(t) = {\rm Mod}(R_T(t)). \]
These two quantities are functionally related,  and the following properties are not difficult to establish.
\begin{lemma}  For $t>1$,  $\gamma_n(t)=2^{n-1}\tau_n(t^2-1)$,  both $\gamma_n$ and $\tau_n$ are continuous, strictly monotone and 
\begin{eqnarray*}  
\lim_{t\searrow 1} \gamma_n(t) = +\infty, && \lim_{t\nearrow \infty} \tau_n(t) = 0\\
  \lim_{t\searrow 0} \tau_n(t) = +\infty, &&  \lim_{t\nearrow \infty} \tau_n(t) = +\infty .
\end{eqnarray*}
\end{lemma}
In two-dimensions these quantities can be explicitly written down in terms of elliptic integrals,  but there are no such formulas known in higher dimensions.  

As will become apparent in a moment,  it is  necessary to get fairly good estimates of these functions at these extreme values as it is from these that equicontinuity results can be deduced,  although there are other ways of course.  Theorem \ref{Rcont} establishes the continuity of these functions. Here are some important estimates on these moduli due to Gehring, \cite{Gehring1}.  They are asymptotically sharp as $t\to\infty$,  but slightly better,  if somewhat more complicated, estimates are known \cite{Vu}.  The number $\lambda_n$ below is known as the Gr\"otzsch ring constant.  The value of $\lambda_n$ is unknown in any dimension other than two,  however we do know $\lambda_{n}^{1/n}\to e$ as $n\to\infty$.
\begin{theorem}  For each $n\geq 2$ there is a constant $\lambda_n\in [4, 2e^{n-1})$,  $\lambda_2=4$,  such that
\begin{eqnarray}
\label{RGest} \omega_{n-1}\big(\log(\lambda_n \, t )\big)^{1-n} \leq &\gamma_n(t)& \leq \omega_{n-1}\big(\log(t) \big)^{1-n} \\
\label{RTest}  \omega_{n-1}\big(\log(\lambda_{n}^{2} \, t )\big)^{1-n} \leq &\tau_n(t-1)& \leq \omega_{n-1}\big(\log(t) \big)^{1-n} .
\end{eqnarray}
\end{theorem}   
What is most important about the Teichm\"uller and Gr\"otzsch rings are the extremal properties.  These are proved by a higher-dimensional generalisation,  due to Gehring \cite{Gehring3},  of the classical technique of symmetrisation  in the complex plane.  This was  first used by Teichm\"uller for these sorts of applications.  

Given $x_0\in {\Bbb R}^n$,  the spherical symmetrisation\index{symmetrisation} $E^v$ of $E$ in direction $v$ is defined as follows;  for $r\in[0,\infty]$,  $E^+\cap {\Bbb S}^{n-1}(x_0,r)\neq \emptyset$ if and only if $E\cap {\Bbb S}^{n-1}(x_0,r)\neq \emptyset$ and then $E^v\cap {\Bbb S}^{n-1}(x_0,r)$ is defined to be the closed spherical cap centred on $x_0+rv$ with the same $(n-1)$-spherical measure as $E\cap {\Bbb S}^{n-1}(x_0,r)$.   Thus if $E$ is connected, then so is $E^v$ and $E^v$ is rotationally symmetric about the ray $x_0+r v$,  $r>0$.  

We symmetrise a ring consisting of two components $E$ and $F$ by symmetrising $E$ in the direction $v$ to get $E^*$  and $F$ in the direction $-v$ to get $F^*$. Then $E^*$ and $F^*$ are disjoint as the spherical measure $(E\cup F)\cap {\Bbb S}^{n-1}(x_0,r)$ is strictly smaller than ${\Bbb S}^{n-1}(x_0,r)$ and so $R^*=R(E^*,F^*)$  is again a ring.  Then we have the following very useful theorem:
\begin{theorem}  Let $x_0\in{\Bbb R}^n$ and $v\in {\Bbb S}^{n-1}$ and let $R=R(E,F)$ be a ring and $R^*=R(E^*,F^{*})$ be its symmetrisation.  Then
\[{\rm Mod}(R^*) \leq {\rm Mod}(R).\]
\end{theorem} 
Next,  a symmetrised ring contains a M\"obius image of a Teichm\"uller ring.  Namely,  the two line segments $E^*\cap L_v$ and $F^*\cap L_v$ in the ray $L_v=\{tv:t\in{\Bbb R}\}$.  Note that  only one of which may be of infinite length,  however they might both be finite.  In the latter case a M\"obis transformation can be used to ensure one of the components is unbounded.  This leads to the following important extremal properties of Teichm\"uller rings.
\begin{theorem}\label{ETR}  Let $R(E,F)$ be a ring with $a,b\in E$ and $c,\infty\in F$.  Then
\[ {\rm Mod}(R) \geq {\rm Mod}\; R_T\Big( \frac{|a-c|}{|a-b|}\Big).\] 
\end{theorem}
By conformal invariance we obtain the corollary
\begin{corollary}   Let $R(E,F)$ be a ring with $a,b\in E$ and $c,d\in F$.  Then
\[ {\rm Mod}(R) \geq {\rm Mod}\big(R_T\big(|a,b,c,d|)\big),  \] 
where the cross ratio\index{cross ratio} is defined by
\begin{equation}
|a,b,c,d|=\frac{|a-c||b-d|}{|a-b||c-d|} .
\end{equation}
\end{corollary}

\subsection{H\"older continuity}
From the extremity of the Gr\"otzsch and Teichm\"uller rings and estimates on their modulus we obtain modulus of continuity estimates for quasiconformal mappings.  Ultimately these give H\"older continuity estimates once we estimate a certain distortion function which we now describe.  Let $f:\Omega\subset {\Bbb R}^n\to{\Bbb R}^n$ be $K$-quasiconformal,  $x\in \Omega$ and put $r=d(x,\partial\Omega)$.  For all $y\in \Omega$ with $|x-y|<r$,  the ring $R={\Bbb B}^n(x,r)\setminus [x,y]$ lies in $\Omega$ and is conformally equivalent to the Gr\"otzsch ring $\gamma_n(r/|x-y|)$ by an obvious inversion.  Next $f(R)$ is a ring with one finite component containing $f(x)$ and $f(y)$ and the other unbounded.  By Theorem \ref{ETR},  the extremailty of the Teichm\"uller ring,  we have 
\[ \tau_n\Big(\frac{d(f(x),\partial\Omega)}{|f(x)-f(y)|} \Big) \leq K \gamma_n\big(\frac{d(x,\partial\Omega)}{|x-y|}\big) \]
Thus 
\begin{equation}\label{2.18}
\frac{|f(x)-f(y)|}{d(f(x),\partial\Omega)}  \leq \varphi_{n,K}\big(\frac{|x-y|}{d(x,\partial\Omega)} \big)
\end{equation}
where $\varphi_{n,K}$ is the distortion function\index{distortion function}
\begin{equation}
\varphi_{n,K}(t) = \tau_{n}^{-1}(K\gamma_n(1/t)).
\end{equation}
Because of (\ref{2.18}) and its various generalisations to special situations,  the distortion function is much studied.  Gehring \cite{Gehring1} showed that there was a constant $C_{n,K}$ such that 
\begin{equation}\label{dsitf}
\varphi_{n,K}(t) \leq C_{n,K} \; t^{1/K} 
\end{equation}
whenever $t<1/2$ but there are much more refined estimates now, see \cite{Vu}.  Combining both (\ref{2.18}) and (\ref{dsitf}) gives H\"older continuity,  and in fact equicontinuity since the constants do not depend on the map in question,  but only their distortion.
  
\begin{theorem} \label{holder} Let $f:\Omega\to\Omega'$ be a homeomorphism such that
\begin{equation}\label{Kde}
{\rm Mod}(f(R)) \leq K {\rm Mod}(R), 
\end{equation}
for all rings $R\subset \Omega$.  Then for all $y< d(x,\partial\Omega)/2$ we have the modulus of continuity estimate
\begin{equation}
\frac{|f(x)-f(y)|}{d(f(x),\partial\Omega')} \leq C_{n,K} \Big(\frac{ |x-y|}{d(x,\partial \Omega)} \Big)^{1/K}.
\end{equation}
\end{theorem} 
In two-dimensions everything can be made rather more precise.  For instance the following is well known.
\begin{theorem}  Let $f:{\Bbb D}\to{\Bbb D}$,  $f(0)=0$ be a $K$-quasiconformal mapping of the unit disk into itself.  Then 
\[ |f(z)| \leq 4^{1-1/K} |z|^{1/K}. \]
\end{theorem}
This theorem has nice asymptotics as $K\to 1$ recovering the classical Schwarz inequality.  Further,  the $K$-quasiconformal map $z\mapsto z|z|^{1-1/K}$ shows the H\"older exponent to be optimal as well.  

\subsection{Mori distortion theorem}
Actually,  if one follows the ideas above and estimates on a larger scale one achieves an important result of Mori \cite{Mori, Vu}
\begin{theorem}  There is a constant $C_{n,K}$ such that if f $f:{\Bbb R}^n\to{\Bbb R}^n$ with $f(0)=0$ is $K$-quasiconformal with respect to rings,  then
\begin{equation}
|f(x)| \leq C_{n,K}\;  |f(y)|
\end{equation}
whenever $|x|=|y|$.  We have the estimate
\[ C_{n,K} \leq [\gamma_{n}^{-1}(\gamma_{n}(\sqrt{2})/K)]^2 \]
where $\gamma_n$  is the Gr\"otzsch ring modulus.
\end{theorem}
When the normalisation $f(0)=0$ is removed we have
\begin{corollary}  If   $f:{\Bbb R}^n\to{\Bbb R}^n$ is $K$-quasiconformal,  then
\begin{equation}
|f(x)-f(z)| \leq C_{n,K}\;  |f(y)-f(z)|
\end{equation}
whenever $|x-z|=|y-z|$.
\end{corollary}
Mori's result is one of a class of results in the distortion theory of the geometry of mappings.  Many other such can be found in Vuorinen's book \cite{Vu},  including higher-dimensional versions of the Schwarz lemma and so forth for quasiregular mappings.  Two further interesting results for quasiconformal mappings measure the distortion of the cross ratio of the points $x,y,z,\infty$.  These in effect lead to the notion of quasisymmetry\index{quasisymmetry} and when M\"obius invariance is used to normalise away the behaviour at $\infty$ we get the notion of quasim\"obius mappings\index{quasim\"obius}.  The ideas are not particularly difficult and follow in much the same way as the distortion estimate of (\ref{2.18}).

\begin{theorem}[local quasisymmetry] \label{qsloc} For each $K\geq 1$ and $s\in (0,1)$, there is a strictly increasing function $\eta_{s,K}:[0,\infty)\to[0,\infty)$ with $\eta(0)=0$ with the following properties.

If $x,y,z\in \overline{{\Bbb B}^n(0,s)}$ with $x\neq z$,  then
\begin{equation}
\frac{|f(x)-f(y)|}{|f(x)-f(z)|} \leq \eta_{s,K}\Big(\frac{|x-y|}{|x-z|} \Big)
\end{equation}
for every $K$-quasiconformal $f:{\Bbb B}^n\to{\Bbb R}^n$.
\end{theorem}
Explicit (but a little complicated) formulas are easily obtained for the function $\eta_{s,K}$ in terms of the Gr\"otzsch and Teichm\"uller functions.  Notice that from (\ref{qsloc}), and the obvious fact that $\eta_{s,K}(1) \geq \eta_{t,K}(1)$ if $s\leq t$, one immediately obtains the bound on the linear distortion
\[ H(x,f) \leq \eta_{s,K}(1).\]
After a rescaling argument, we also obtain a global version of this.
\begin{theorem}[global quasisymmetry]  For each $K\geq 1$ there is a strictly increasing function $\eta_K:[0,\infty)\to[0,\infty)$ with $\eta(0)=0$ such that 
\begin{equation}
\frac{|f(x)-f(y)|}{|f(x)-f(z)|} \leq \eta_K\Big(\frac{|x-y|}{|x-z|} \Big)
\end{equation}
for every $K$-quasiconformal $f:{\Bbb R}^n\to{\Bbb R}^n$.
\end{theorem}

\section{Compactness} 
An equally important aspect of quasiconformal mappings are their compactness properties.  Usually these are couched in terms of normal family type results.  Recall a family of mappings 
\[ {\cal F}= \{f:\Omega\subset{\Bbb R}^n\to{\Bbb R}^n\} \]
is said to be normal\index{normal family} if every sequence $\{f_n\}_{n=1}^{\infty}\subset{\cal F}$ contains a subsequence which converges uniformly on compact subsets of $\Omega$.  The modulus of continuity estimates of Theorem \ref{holder} guarantee the equicontinuity,  and therefore via the Arzela-Ascoli theorem,  the normality of any suitably normalised family of $K$-quasiconformal mappings.  Not only that of course,  the bilipschitz estimate on the distortion of moduli also shows the family of inverses (restricted to a suitable domain of common definition) is also normal. These observations quickly lead to compactness
results. The most elementary of these is the following.

\begin{theorem}\label{normF} Let $\Omega\subset{\Bbb R}^n$ and $x_0,y_0\in \Omega$.  Then the family
\[{\cal F}_K=\{f:\Omega\to{\Bbb R}^n, \; f(x_0)=0, f(y_0)=1, \; \mbox{ and $f$ is $K$-quasiconformal} \} \]
is a normal family.
\end{theorem}
Obviously something is necessary here as the family of conformal mappings of ${\Bbb R}^n$,   $\{x\mapsto n x\}_{n=1}^{\infty}$  is not normal.  If one wants to add the point at $\infty$ to the discussion and  consider families of $K$-quasiconformal maps defined on the Riemann sphere $\hat {{\Bbb R}}^n$ normalisation at three points is all that is required to guarantee normality.  Next,  as convergence in given by Theorem \ref{normF} the issue is whether the limit map is actually quasiconformal.  The next theorem establishes this.

\subsection{Limits of quasiconformal mappings}
\begin{theorem}\label{cthm2} Fix $K$.  Let $f_j:\Omega \to \Omega_j$ be a sequence of $K$--quasiconformal mappings converging pointwise to
$f:\Omega
\to
{\Bbb R}^n$.  Then one of the following occurs.
\begin{itemize}
\item  $f$ is a $K$--quasiconformal embedding and the convergence is uniform on compact subsets.
\item  $f(\Omega)$ is a doubleton with one value attained only once.
\item  $f$ is constant.
\end{itemize}
\end{theorem}
When $\Omega = \hat {{\Bbb R}}^n$,  and with the obvious interpretation of continuity at infinity and so forth,  we have the following convergence properties of quasiconformal mappings of the Riemann sphere.

\begin{theorem}\label{cthm3}Fix $K$.  Let $f_j:\hat {{\Bbb R}}^n \to \hat {{\Bbb R}}^n$ be a sequence of $K$--quasiconformal mappings.  Then there is a subsequence $\{f_{j_k}\}_{k=1}^{\infty}$ such that  one of the following occurs;
\begin{itemize}
\item  there is a $K$ quasiconformal homeomorphism $f:\hat {{\Bbb R}}^n\to\hat {{\Bbb R}}^n$  and both $f_{j_k}\to f$ and $f_{j_k}^{-1}\to f^{-1}$ uniformly on $\hat {{\Bbb R}}^n$,  or
\item  there are constant $x_0,y_0\in \hat {{\Bbb R}}^n$,  possibly $x_0=y_0$,  such that
\begin{eqnarray*}
 f_{j_k} \to x_0 &&  \mbox{locally uniformly in $\hat {{\Bbb R}}^n\setminus\{y_0\}$}, {\rm and}\\
 f_{j_k}^{-1} \to y_0 &&  \mbox{locally uniformly in $\hat {{\Bbb R}}^n\setminus\{x_0\}$} .
 \end{eqnarray*}
\end{itemize}
\end{theorem}

Basically it is true that any sufficiently normalised family of quasiconformal mappings forms a normal family. 
However,  there is a far reaching generalisation of these sorts of results.  It is Rickman's version of Montel's Theorem which we
discuss a bit below. 

\medskip

There are also more general results about the normal family properties of families of mappings with finite distortion.  Typically very little can be said but if,  for instance,  the distortion function 
\begin{equation}\label{distortionfunction} K(x,f) = \frac{|Df(x)|^n}{J(x,f)}
\end{equation}
has strong integrability properties such as being exponentially integrable,  then there are very similar results to those of Theorems \ref{normF}, \ref{cthm2} and \ref{cthm3} available,  \cite{IM}.

\section{The analytic definition of a quasiconformal mapping}

Examining the proof of Theorem \ref{mci},  it becomes apparent that we should get the Lipschitz estimate of (\ref{Moddist}) if we were to have the pointwise estimate between the differential matrix $Df$ and its determinant.
\begin{equation}\label{Keq}
|Df(x)|^n \leq K \; J(x,f)
\end{equation}
which of course is close to (\ref{confdef}).  We would need this for both $f$ and its inverse of course,  but at least where the differential is nonsingular if we write out the eigenvalues of the symmetric positive definite matrix $Df^t(x)Df(x)$ as $\lambda_1\leq \lambda_2 \leq \cdots \leq \lambda_n$,  then the inequality (\ref{Keq}) reads as 
\[ \lambda_{n}^{n} \leq K^2\;  (\lambda_1 \times \lambda_2 \times\cdots\times \lambda_n)  \]
and this certainly implies $\lambda_n \leq K^2 \lambda_1$ and hence
\[   \lambda_1 \times \lambda_2 \times\cdots\times \lambda_n \leq  \lambda_1 \lambda_{n}^{n-1}  \leq K^{2(n-1)} \lambda_{1}^{n} .\]
Therefore writing $g=f^{-1}:\Omega'\to\Omega$ we would have
\begin{equation}\label{Keq2}
|Dg(y)|^n \leq K^{n-1} \; J(y,g), \hskip10pt y=f(x)
\end{equation}
so $g$ will also have a Lipschitz estimate,  thus giving the bilipschitz estimate we want - albeit with different constants $K$.   One of course needs some sort of Sobolev regularity to make all this work,  and that leads us to the analytic definition of quasiconformal mappings. \index{quasiconformal!analytic definition}
Let $f:\Omega \to \Omega^\prime$ be a homeomorphism belonging to the Sobolev class $W^{1,n}_{loc}(\Omega)$ of functions whose first derivatives are locally $L^n$-integrable.  Then $f$ is $K$--quasiconformal if there exists a $K$, $1 \leq K \leq
\infty$,  such that
\begin{equation} \label{adef}
|Df(x)|^n \leq K \; J(x,f), \hskip20pt \mbox{ at almost every point $x\in \Omega$.}
\end{equation}
 We again need to point out that the constant $K$ here is not the same as that for rings.  Further,  it is not in general true that the composition of $W^{1,n}$ mappings is again of Sobolev class $W^{1,n}$,  nor is it true that the inverse of a $W^{1,n}$ homeomorphism is $W^{1,n}$,  so the fact that the composition of quasiconformal mappings and inverse the inverse of a quasiconformal mapping are again quasiconformal,  discussed in Theorem \ref{composition},   are now not nearly so direct.
 
 \medskip

There are advantages however;  considering this definition,  one sees that the hypothesis that $f$ is a homeomorphism to be largely redundant.  We therefore say that $f$ is $K$-{\em quasireguar}\index{quasiregular} if $f$ satisfies (\ref{adef}) and has the appropriate $W^{1,n}_{loc}$ Sobolev regularity.  In fact the hypothesis that $f\in W^{1,n}_{loc}(\Omega)$ ensures that the Jacobian determinant of $f$ is a locally integrable function and gives one a chance of establishing such things as the change of variable formula and so forth.  From this purely
analytic definition of quasiconformal mappings,  Reshetnyak was able to establish important topological properties.\index{Reshetnyak's theorem}

\begin{theorem}\label{RST}{\rm (Reshetnyak)}  A quasiregular mapping $f:\Omega\subset{\Bbb R}^n\to{\Bbb R}^n$ is open and discrete 
\end{theorem}
With this we now recall Chernavskii's theorem \cite{Chern}.
\begin{theorem} \label{Top} Let  $B_f$ denote the branch set of a quasiregular mapping $f:\Omega\to{\Bbb R}^n$,  that is
\begin{equation}
B_f = \{x\in\Omega:\mbox{$f$ is not locally injective at $x$}\}.
\end{equation}
Then  the topological dimension of both $B_f$ and $f(B_f)$ is less than or equal to $n-2$.
\end{theorem}
These two results,  Theorems \ref{RST} and \ref{Top},  together with fairly general topological degree theory and covering properties of branched open mappings established by V\"ais\"al\"a \cite{Vaisala1} and others give  various quite general path lifting properties of these mappings  \cite{MRV1,MRV2,MRV3,R}, and the
well-known result of Poletsky   \cite{Poletsky}.  With these properties at hand one can study the deformations of curve families in the more general setting of quasiregular mappings.  The distortion bounds on the modulus  enable the geometric methods of modulus to be used to great effect to build a theory analogous to that of analytic functions in the complex plane.  

\medskip

There is a considerable body of research building around these topological questions for mappings of finite distortion.  The questions become deep and subtle and beyond the scope of this chapter,  but the interested reader can consult \cite{IM} and work forward to the many interesting current research directions.

\medskip

There are a few other consequences of the analytic definition that need to be recounted.  These are key features for the analytic theory of these mappings showing sets of zero-measure are preserved.
\begin{theorem}[Condition $N$ and $N^{-1}$]
Let $f:\Omega\to{\Bbb R}^n$ be non-constant quasiregular mapping.  
\begin{itemize}
\item If $A$ has measure $0$,  $|A|=0$,  then  $|f(A\cap \Omega)|=0$.
\item If $|B|=0$,  then $|f^{-1}(B)|=0$.
\item $J(x,f)>0$ almost everywhere in $\Omega$.
\item $|B_f|=0$,  and hence $|f(B_f)|=0$.
\end{itemize}
\end{theorem}

\section{The Liouville Theorem} 

In 1850, the celebrated French mathematician Joseph Liouville added a 
short note to a new edition of Gaspard Monge's classic work 
{\em Application de l'Analyse \`a la G\'eometrie\/}, whose publication 
Liouville was overseeing.
The note was prompted by a series of three letters that Liouville had 
received in  1845 and 46 from the renowned British physicist 
William Thomson.  Thomson, better known today as Lord Kelvin, had studied in Paris under 
Liouville's  in the mid-1840s, so these two 
giants of nineteenth century science were well acquainted. 

In his letters, Thomson asked Liouville a number of questions concerning
inversions in spheres, questions that had arisen in conjunction with 
Thomson's research in electrostatics, in particular, with the so-called 
principle of electrical images 
(we point out that the reflection in the unit sphere ${\Bbb S}^2$ of ${\Bbb R}^3$
is often referred to in physics as the ``Kelvin transform.'')
More about the interesting relationship between Thomson and 
Liouville can be found in Jesper L\"utzen's 
magnificient biography of Liouville .

The substance of 
Liouville's note is conveyed by the following remarkable assertion: \index{Liouville theorem}

\begin{theorem}[Liouville's theorem]
If $\Omega$ is a domain in ${\Bbb R}^n$,  $n \geq 3$, then any conformal mapping 
$f:\Omega \to {\Bbb R}^n$ is the restriction to $\Omega$ of a M\"obius transformation of $\hat {{\Bbb R}}^n$. 
\end{theorem}

Exactly what is meant by a conformal mapping is a modern day issue around the regularity theory of solutions to PDEs such as the Cauchy-Riemann system below at (\ref{CRE}).  But in Liouville's time he certainly understood such mappings to be many times differentiable and following his motivation for writing the article, Liouville couched 
his discussion in the language of differential forms rather than mappings. 
As a consequence, his original formulation  bears little resemblance
to the theorem above, although the relationship between the two 
formulations is quite clear via
differential geometry. However Liouville's title,  ``Extension au cas de trois dimensions de la question du 
trac\'e g\'eographique''  gives no hint whatsoever as to the results. 
 
It was only later that Liouville published his theorem in a form 
approximating the statement of it that we have given.

As we hinted,  the proof which Liouville 
outlined for his theorem makes use of certain implicit smoothness hypotheses which when unwound gives the added assumption that $f$ is mapping of class $C^3(\Omega)$,   that is three times continuously differentiable,  or better.  In higher dimensions the group {\bf M\"ob}($n$) of all M\"obius transformations\index{M\"obius transformations} of
$\hat {{\Bbb R}}^n$  consists of all finite compositions of reflections in spheres and hyperplanes.  It is easy to see that these mappings provide 
examples of conformal transformations.  They are of course all $C^{\infty}(\Omega)$.  

The smoothness assumption in Liouville's theorem is not optimal.  One would like to relax the injectivity
assumption to allow the possibility of branching and also to relax the differentiability assumption
as much as possible.  The natural setting for Liouville's theorem is a statement about the regularity of solutions to the {\em Cauchy-Riemann} system\index{Cauchy-Riemann system} 
\begin{equation}\label{CRE}
D^tf(x)\; Df(x) = J(x,f)^{2/n} {\rm Id}, \hskip20pt \mbox{almost everywhere in $\Omega$},
\end{equation} 
where ${\rm Id}$ is the $n\times n$ identity matrix.   If $f:\Omega\to{\Bbb R}^n$ is a $1$-quasiconformal mapping,  then pointwise almost everywhere we must have the positive semidefinite matrix $Df^t(x)Df(x)$ having a single eigenvalue with multiplicity $n$.  Thus either $Df(x)=0_{n\times n}$ and $J(x,f)=0$ or,  as a little linear algebra will reveal,  $Df(x)$ is a scalar multiple of an orthogonal transformation.   In particular,  with the analytic definition of quasiconformality we see that a $1$-quasiconformal mapping $f$ is a $W^{1,n}_{loc}(\Omega)$ solution to the equation (\ref{CRE}).

With this formulation we have the following very strong version of the Liouville theorem established using the nonlinear Hodge Theory developed in \cite{IM1,I2}. In two-dimensions it is analogous to the classical Looman-Menchoff Theorem. Note especially that there is no longer any assumption of injectivity - it is a consequence of the theorem.

\begin{theorem}\label{LT} Let $\Omega\subset {\Bbb R}^n$, $n\geq 3$ and let $f\in W^{1,p}_{loc}(\Omega)$ be a weak solution to the equation (\ref{CRE}). If $n$ is even,  then any solution with $p\geq n/2$ is the restriction to $\Omega$ of a M\"obius transformation of $\hat {{\Bbb R}}^n$.  If $n$ is odd,  then there is an $\epsilon=\epsilon(n)$ such that any weak solution with $p\geq n-\epsilon$ is the restriction to $\Omega$ of a M\"obius transformation of $\hat {{\Bbb R}}^n$.

\medskip

This is sharp in the following sense.  In all dimensions $n\geq 2$ and all $p<n/2$,  there is a weak solution to (\ref{CRE}) in $W^{1,p}_{loc}(\Omega)$  which is not in $W^{1,n}_{loc}(\Omega)$ and so in particular is not a M\"obius transformation.
\end{theorem}
 
The discrepancy here between what is known in odd dimensions and even dimensions is one of the central unsolved problems in the theory.  Further,  although the results are very sharp in even dimensions,  there remains the possibility of improvement.  For instance it might be that the Liouville theorem remains true in ${\Bbb R}^n$,  $n\geq 3$,  for weak $W^{1,1}_{loc}(\Omega)$ solutions which are continuous.

\bigskip

We draw attention to one significant corollary of Liouville's theorem:
the only subdomains $\Omega$ in ${\Bbb R}^n$ with $n\geq 3$ that are conformally 
equivalent to the unit ball ${\Bbb B}^n$ are Euclidean balls and half-spaces. 
This stands in stark contrast to the marvelous discovery by Riemann, announced in 1851 
a year after Liouville's note was published:
any simply connected proper subdomain $\Omega\subset {\Bbb C}$ of the complex 
plane is conformally equivalent to the unit disk ${\Bbb D}$. 

 \bigskip
 
 From the formulation of Liouville's theorem in Theorem \ref{LT} we are naturally led to the basic
connections between quasiregular mappings and non-linear PDEs  through  the {\em Beltrami
system}. \index{Beltrami system} Let  $S(n)$ denote the space of symmetric  positive definite $n \times n$
matrices of determinant equal to $1$.  Geometrically $S(n)$ is  a non positively curved complete symmetric space. 

Given $\Omega$ a subdomain of ${\Bbb R}^n$ and $G:\Omega \rightarrow S(n)$ a
bounded measurable mapping we define the Beltrami equation as 
\begin{equation}\label{BE}
D^tf(x)\; Df(x) = J(x,f)^{2/n} G(x), \hskip10pt\mbox{almost every $x\in\Omega$}
\end{equation}

 To each non-constant quasiregular mapping, 
there corresponds a unique (tautological) Beltrami equation and we refer to $G$ as the {\em distortion tensor}\index{distortion tensor}  of the mapping $f$. 

  A key approach to the modern theory is to examine properties and obtain geometric information about quasiregular mappings (and more general mappings of
finite distortion) when they are viewed as solutions to this and related PDEs. These equations are studied from many points of view, as the Euler--Lagrange equations for the {\em absolute} minima of  variational integrals,  at the level of
differential forms using exterior algebra and also as equations relating the Dirac operators of conformal and spin
geometry,  see \cite{IM} for results in these directions.

\section{Gehring's Higher Integrability}\index{higher integrability}

In a remarkable paper in 1973,  F.W. Gehring established that the Jacobian determinant of a
$K$-quasiconformal mapping is integrable above the natural exponent.  That is the assumption $f\in
W^{1,n}_{loc}(\Omega,{\Bbb R}^n)$  together with the bound on distortion implies that $f\in
W^{1,n+\epsilon}_{loc}(\Omega,{\Bbb R}^n)$ for some $\epsilon$ depending on $n$ and $K$. Gehring gave explicit estimates on $\epsilon$. While this result was already known in the plane due to the work of Bojarski
\cite{B}, and perhaps anticipated in higher dimensions, it is impossible to overstate how important this result 
has proven to be in the theory of quasiconformal mappings and more generally Sobolev spaces and in the
theory of non-linear PDEs.  The techniques developed to solve this problem,  for instance the
well--known reverse H\"older inequalities are still one of the main tools used in several disciplines,
including non-linear potential theory, non-linear elasticity, PDEs and harmonic analysis. 

We state the following version of Gehring's result as proved in \cite{IM} which also gives the result  for quasiregular mappings.

\begin{theorem}  Let $f:\Omega\to{\Bbb R}^n$ be a mapping of Sobolev class $W^{1,q}_{loc}(\Omega)$ satisfying the differential inequality
\begin{equation}\label{difineq}
|Df(x)|^n \leq K\; J(x,f).
\end{equation}  
Then there are $\epsilon_{K}{*},\epsilon_K>0$ such that if $q>n-\epsilon_{K}{*}$,  then $f\in W^{1,p}_{loc}(\Omega)$ for all $p<n+\epsilon_K$.
\end{theorem}
As an immediate corollary we have the following:
\begin{theorem}\label{GHI}  Let $f:\Omega\to{\Bbb R}^n$ be a $K$ quasiconformal mapping.  Then there is $p_{_K}>n$ such that  $f\in W^{1,p_{_K}}_{loc}(\Omega)$.
\end{theorem}
The higher-dimensional integrability conjecture here would assert that if $f$ satisfies (\ref{difineq}) and lies in $W^{1,q}_{loc}(\Omega)$ for some $q>nK/(K+1)$, then $f$ actually lies in the Sobolev space $W^{1,p}_{loc}(\Omega)$ for all $p<nK/(K-1)$.  In two-dimensions this conjecture was proven by K. Astala \cite{Astala}.  In even dimensions rather more is known and the numbers $\epsilon_{K}{*}$ and $\epsilon_K>0$ can be related to the $p$-norms of certain singular integral operators which can be estimated.  Indeed the conjecture would follow from the current conjectural identification of these norms.  In odd dimensions rather less is known.  In any case,  Theorem \ref{GHI} yields the following\index{reverse H\"older inequality}
\begin{corollary}[Reverse H\"older inequality]  Let $f:\Omega\subset {\Bbb R}^n \to \Omega'$ be a $K$-quasiconformal mapping.  Then there is $p=p(n,K)>1$ and $C=C(n,K)$ such that
\begin{equation}\label{RHI}
\Big(\frac{1}{|Q|}\, \iint_Q \; J(x,f)^p \; dx \Big)^{1/p} \leq \frac{C}{|Q|}\, \iint_Q \; J(x,f) \; dx
\end{equation}
for all cubes $Q$ such that $2Q\subset\Omega$.
\end{corollary}
Actually,  our presentation here is a bit back to front as it is via the reverse H\"older inequality at (\ref{RHI}) that the higher integrability Theorem \ref{GHI} was first established.
The restriction to cubes $Q$ so that $2Q\subset\Omega$ is necessary but can be removed under assumptions about the regularity of $\Omega'=f(\Omega)$ - namely that it should be a John domain.  There is another connection here to the nonlinear potential theory as the estimate shows the Jacobian $J(x,f)$ to be an $A_\infty$ Muckenhoupt weight on the cubes $Q$.

\medskip

Another interesting unsolved problem concerns the question of when a positive function can be the Jacobian of a quasiconformal mapping.  Obviously the results above impose restrictions on such a function.  Another is given by Reimann's result:\index{Reimann's theorem}
\begin{theorem}[Reimann's theorem]  Let $f:{\Bbb B}^n\to{\Bbb R}^n$ be $K$-quasiconformal.  Then $\log J(x,f)$ is a function of bounded mean oscillation.
\end{theorem}

\section{Further Stability and Rigidity Phenomena}

Along with the Liouville theorem there are other interesting phenomena which occur only in higher dimensions, $n\geq 3$.  For
instance consider the following local to global homeomorphism property.  In
1938 Lavrentiev \cite{Lav} asserted that a locally homeomorphic quasiconformal mapping
${\Bbb R}^3 \to {\Bbb R}^3$ is a global homeomorphism onto.  This assertion was
proved correct by Zorich \cite{Zorich} in all dimensions $n\geq 3$.

\begin{theorem}
Let $f:{\Bbb R}^n \to {\Bbb R}^n$ be a locally homeomorphic quasiregular mapping.  If $n\geq 3$, then $f$ is a globally
injective quasiconformal mapping onto ${\Bbb R}^n$.
\end{theorem}

The condition $n\geq 3$ is essential as the exponential mapping $e^z$ in the plane demonstrates.  Zorich's theorem\index{Zorich's theorem} was
generalised by Martio--Rickman--V\"ais\"al\"a  in the following way (an
earlier result of John  proved the same result for locally bilipschitz mappings).

\begin{theorem}
There is a positive constant $r=r(n,K)$ with the following property.  If  $f:{\Bbb B}^n \to {\Bbb R}^n$ is a locally injective
$K$--quasiregular mapping,  then  $f|{\Bbb B}^n(0,r)$ is injective.
\end{theorem}

The number $r(n,K)$ in the above theorem is called the injectivity radius.  Zorich's result clearly follows from this result
by scaling.

\medskip

There are also interesting local to global injectivity results for quasiregular mappings between Riemannian manifolds.  In
this vein the following result of Gromov is perhaps best known \cite{Gromov}.

\begin{theorem}  If  $f:M\to N$ is a locally homeomorphic quasiregular mapping of a complete Riemannian
$n$-manifold  $M$  of finite volume into a simply connected Riemannian manifold $N$  with
$n\geq 3$,  then  $f$ is injective and $N\backslash f(M)$ is of Hausdorff dimension
zero.
\end{theorem}

These results are very well presented in \cite{R}.  There are also stability results of a different nature. 
These are based on the compactness properties of mappings of finite distortion and the Liouville theorem. Roughly speaking
one can show that in all dimensions as $K\rightarrow 1$,  $K$-quasiregular mappings are
uniformly well approximated by conformal mappings.  Liouville's theorem implies in dimension $n\geq
3$ that conformal mappings are M\"obius transformations.  Thus in dimension $n\geq 3$ for sufficiently
small $K$ we obtain local injectivity by virtue of the uniform approximation by a globally injective
mapping,  see \cite{Reshetnyak1} and for an interesting
application \cite{MS}. For instance one has
\begin{theorem}  For each $n\geq 3$ there is a constant $\delta(n,\epsilon)$ with the following
properties
\begin{itemize}
\item $\delta(n,\epsilon) \rightarrow 0$ as $\epsilon \rightarrow 0^+$
\item If  $f:{\Bbb B}^n \rightarrow {\Bbb R}^n$
is a $K$--quasiregular mapping with $K\leq 1+\epsilon$,  then there is a M\"obius mapping  $\phi:{\Bbb B}^n
\rightarrow {\Bbb R}^n$ such that \begin{equation}
\sup_{x\in {\Bbb B}} |(\phi^{-1}\circ f)(x)-x| < \delta(n,\epsilon)
\end{equation}
\end{itemize}
\end{theorem}

Also we mention the following connection between distortion and local injectivity.  Again,  this is a higher-dimensional phenomena.  The map $f:z\mapsto z^2$ is a $1$-quasiregular map of ${\Bbb C}$ with $B_f=\{0\}\neq \emptyset$.  Contrast this with the following theorem.

\begin{theorem}\label{MartioC}  There is a constant $K_0> 1$ with the following property. Let $f:\Omega\to{\Bbb R}^n$,  $n\geq 3$, be $K$-quasiregular.  If $K<K_0$,  then $B_f=\emptyset$. That is $f$ is locally injective.
\end{theorem}
The number $K_0$ depends on the particular definition of the distortion $K$.  But with the geometric and analytic definitions the number is expected to be equal to $2$ (known as the Martio conjecture)\index{Martio conjecture}.  The best known bound is due to Rajala \cite{Rajala} and is only very slightly bigger than $1$, but it is explicit and not derived from a compactness argument.  It would be a major advance to establish the sharp result here.

\medskip

Among other consequences it is known that if
the distortion tensor of a quasiregular mapping is close to continuous in the space of functions of bounded mean oscillation (BMO),  or in particular continuous,  then local injectivity
follows.  Closely related results can be found in \cite{MRSY}.   Again,  many of these sorts of results are based around compactness arguments and do not give effective information.

\medskip

These results explain why we really need to consider measurable conformal
structures in the defining equation (\ref{BE})  as any
degree of regularity of the distortion tensor  forces local injectivity.  Precisely what  regularity is necessary is   a study currently under intense investigation.

\section{Quasiconformal Structures on Manifolds}

In this chapter we will not delve too deeply into the theory of quasiconformal mappings on manifolds.  Of course the local theory, regularity, and compactness results pretty much follow from the Euclidean theory,  but there are some quite subtle and interesting aspects to the theory that warrant deeper investigation.  These investigations are far from complete at present.

The starting point for questions concerning quasiconformal mappings and structures on manifolds is 
 Sullivan's uniformisation theorem \index{Sullivan's uniformisation theorem}\index{uniformisation theorem!Sullivan} which tells us that, apart from dimension $n=4$, every topological manifold admits a unique quasiconformal structure\index{quasiconformal structure}
(quasiconformal coordinates).    This theorem is quite remarkable in that it allows analytical calculation on topological manifolds - manifolds which may not even admit a differentiable structure.  Thus one may seek to calculate topological invariants analytically.  Further,    it is a consequence of uniqueness that two different smooth structures on the same compact manifold are quasiconformally equivalent.  

\medskip

The classical uniformisation theorem in complex analysis states that every surface $F$ admits a
conformal structure.  That is a set of local coordinates
$\{(\varphi_\alpha,U_\alpha)\}_{\alpha\in A}$ with $\bigcup_\alpha U_\alpha = F$ and
$\varphi_\alpha:U_\alpha \hookrightarrow {\Bbb C}$ in which the transition mappings
$ \varphi_{\alpha}\varphi_{\beta}^{-1} : \varphi_{\beta}(U_\alpha \cap U_\beta)\rightarrow {\Bbb C}$
are conformal mappings for all $\alpha$ and $\beta$ between planar subdomains.  Every surface has a
simply connected covering space which inherits this conformal structure.  The monodromy theorem
then implies that this covering space is one of ${\hat{{\Bbb C}}}$, ${\Bbb C}$ or the unit disk ${\Bbb D} = \{z\in
{\Bbb C}:|z|<1\}$.     For every (orientable) surface except ${\hat{{\Bbb C}}}, {\Bbb C}, {\Bbb C}\backslash \{z_0\}$ and the torus, the
universal cover is the unit disk and the group of cover translations is a subgroup of the group of
conformal automorphisms of  ${\Bbb D}$,  that is a group of linear fractional transformations, called
a {\em Fuchsian group}.  

This result is of course one of the most profound results in complex
analysis.  The theory of Fuchsian groups developed by
Poincar\'e  laid the foundations for the study of discrete groups of transformations of
more general spaces and geometries.

Quasiconformal mappings play an essential r\^ole in the study of Fuchsian groups and their orbit
spaces,  {\em Riemann surfaces}.  The theory of Teichm\"uller spaces uses quasiconformal mappings to
study the various conformal structures on a given Riemann surface.  This is amply demonstrated in the contents of this book.

However,  what we want to consider here is the extent to which the uniformisation theorem might be true in
higher dimensions.  Because of the rigidity of conformal mappings in space it is not to be expected
that every $n$-manifold admits a conformal structure. Although Perelman's recent proof of Thurston's geometrisation theorem, building on earlier work of Hamilton,  suggests that this is nearly the case in dimension $3$.

\subsection{The existence of quasiconformal structures}

In general, given any pseudo--group of homeomorphisms of Euclidean space one can define the associated category
of manifolds using the pseudo--group to provide local coordinates.  The two most familiar pseudo--groups are of
course the pseudo--group of homeomorphisms,  giving rise to topological manifolds,  and the
pseudo-group of diffeomorphisms,  giving rise to $C^\infty$ or smooth manifolds.  Other examples of
possible structures would be piecewise linear,  real analytic,  complex analytic and so forth.
One of the fundamental problems of topology has been to determine when a topological manifold
admits a ``nicer'' structure than that given {\em a priori} and how many different sorts of similar structures exist on a particular
manifold.  For instance one might ask: when does a topological manifold admit a smooth structure? Given
a manifold with two potentially different  smooth structures,  are they the same by a smooth change of coordinates? 

Notice that the possibility of admitting a smooth structure is a {\em topological} invariant.  That
is if $M$ and $N$ are homeomorphic and $M$ is smooth,  then $N$ admits a smooth structure,  obtained
by simply declaring that the homeomorphism is a smooth map.  Due to the work of  Moise and others in three dimensions,  the differences between smooth and
topological structures first shows up in dimension $4$. Because of the work of Freedman
\cite{Freedman} and Donaldson \cite{Donaldson} we know that there are plenty
of $4$ manifolds which do not admit any smooth structure and,  quite surprisingly,  topological
manifolds as simple as $4$--dimensional Euclidean space which admit many
different smooth structures.  A similar situation persists in higher
dimensions.

Notice that in order to do calculus or study function theory on a manifold some smoothness
assumptions are necessary on the coordinate charts.   From the geometric point of view  {\em quasiconformal
manifolds} would seem a natural starting point.  We say a manifold  $M$ admits a {\em quasiconformal
structure}  if there is a set of local coordinates
$\{(\varphi_\alpha,U_\alpha)\}_{\alpha\in A}$ with
$\bigcup_\alpha U_\alpha = M$ and $\varphi_\alpha:U_\alpha \hookrightarrow {\Bbb R}^n$ in which the transition mappings
\begin{equation}
\varphi_{\alpha}\varphi_{\beta}^{-1} : \varphi_{\beta}(U_\alpha \cap U_\beta)\rightarrow {\Bbb R}^n
\end{equation} 
are quasiconformal mappings of subdomains of ${\Bbb R}^n$ for all $\alpha$ and $\beta$.  Notice that
there is no assumption on the distortion of the transition charts other than boundedness.  Since quasiconformal mappings of subdomains of ${\Bbb R}^n$
have $L^n$--integrable first derivatives they admit enough structure so as to be able to
define differentiation,  speak of differential forms and exterior derivatives, define a de Rham
type cohomology theory and discuss the index theory of certain differential operators.  We can speak of
conformal and quasiregular mappings between quasiconformal manifolds and study conformal invariants
of such manifolds.  The reader should be aware of the complexity and some nuances of the very definition of the various
Sobolev classes of mappings between manifolds with measurable metric tensors, see Bethuel \cite{Bethuel} and Hajtasz
\cite{Haj}. The circle of ideas concerning the question of regularity of topological manifolds from the point of view of
analysis is interesting and important.

\medskip

There are two principal properties of a pseudo--group of transformations in a given category,  denoted CAT,  (for instance
smooth,  piecewise linear or quasiconformal) to imply that a topological manifold admits such a structure,  and if it does so,
then it is unique.  These are

\medskip

\begin{itemize}
\item  {\em Deformation}.  Two CAT homeomorphisms which are uniformly close in the $C^0$-topology can be deformed one to
the other through CAT homeomorphisms (a suitable relative version of this statement is also necessary).\\
\item {\em Approximation}.  Any homeomorphism $ {\Bbb B}^n \hookrightarrow {\Bbb R}^n$ can be uniformly approximated in the $C^0$-topology by a CAT homeomorphism.
\end{itemize}

In a remarkable piece of work D. Sullivan established the deformation property in all dimensions for the category of
quasiconformal mappings \cite{Su1}  and he also laid the foundations for much of the recent work in geometric topology
in the quasiconformal category,  notably the work of Tukia and V\"ais\"al\"a  \cite{TV}. 
Sullivan also established the approximation property for $n\neq 4$.  The basic tool was a hyperbolic version of the
Edwards--Kirby furling technique of geometric topology.  As a consequence of Sullivan's work we have the following
remarkable result which one might regard as an analogue of the $2$--dimensional uniformisation theorem.

\begin{theorem}
Every topological $n$--manifold,  $n \neq 4$,  admits a unique quasiconformal structure.
\end{theorem}

\subsection{Quasiconformal $4$-manifolds}The revolution  in our understanding of the theory of $4$--manifolds initiated by Donaldson and Freedman has not left the theory of quasiconformal
mappings untouched. As we have discussed above Sullivan's uniformisation theorem implies that every topological
$n$--manifold, $n\neq 4$, admits a quasiconformal structure.  This leaves open the question of what possible structures can
exist on an arbitrary topological $4$--manifold.  Donaldson and Sullivan \index{Donaldson--Sullivan} attacked this
problem in a beautiful paper in 1990
\cite{DS} which heralded many new ideas into the theory of quasiconformal mappings.  Their approach was to take a
quasiconformal $4$-manifold and develop the associated global Yang--Mills \index{Yang--Mills} theory on such a manifold
and thereby produce the same sorts of invariants associated to intersection forms that are used to distinguish the
topological manifolds which admit smooth structures from those that do not.

Recall that in 1982 Freedman gave a complete classification of compact simply connected $4$--manifolds
\cite{Freedman} by establishing the $4$--dimensional $h$--cobordism
theorem.  Thus there is exactly one simply connected topological
$4$--manifold for each given unimodular intersection form.  In 1983 Donaldson  \cite{Donaldson} showed the only
negative definite forms which are realised as the intersection forms
of smooth compact simply connected  $4$--manifolds  are the standard
diagonalisable forms (the hypothesis on simple connectivity was later
removed).  These results then provide a mechanism for finding topological
$4$--manifolds which do not admit smooth structures.  If one could
develop the necessary Yang--Mills theory for manifolds with less
smoothness assumptions,  then one could similarly provide  examples of
topological
$4$--manifolds which do not admit quasiconformal structures.

The development of this theory in the quasiconformal category is highly
non-trivial and significant technical obstructions need to be overcome.  Notice that for
instance quasiconformal manifolds do not admit Riemannian metrics and the smooth construction
depends on splitting the curvature into self--dual and anti--self--dual parts to define the
anti--self-dual moduli space of connections modulo gauge equivalence.  These are the objects from
which the invariants are computed.

In the quasiconformal category,  Donaldson and Sullivan set up some differential geometric
invariants on a quasiconformal manifold $M$ based around the existence of a  measurable conformal
structure.  Principally these were the anti-self dual Yang--Mills equations.  Since the Yang--Mills equations are conformally invariant,  the measurable conformal
structure can be used to define the anti--self--dual connections.  The analysis of these connections
requires the non-linear Hodge theory and the improved regularity properties of quasiregular
mappings.  The fact that quasiconformal mappings preserve the ``correct''  Sobolev spaces plays no
small part in this development. Their main results are as follows.
\begin{theorem}
There are topological $4$--manifolds which do not admit any quasiconformal structure.
\end{theorem}
\begin{theorem}
There are smooth compact (and therefore quasiconformal) $4$--manifolds which are homeomorphic but
not quasiconformally homeomorphic. 
\end{theorem}
As far as we are aware the question of whether there are quasiconformal $4$--manifolds which do not
admit smooth structures remains open.  Also as a consequence of the deformation and approximation
theory properties we discussed above with regard to Sullivan's uniformisation theorem we obtain the
following curiosity in dimension $4$.

\begin{corollary}
There is an embedding $\varphi$ of the $4$--ball,  $\varphi :{\Bbb B}^4  \hookrightarrow {\Bbb R}^4$,  which cannot be approximated uniformly in the spherical metric
by a quasiconformal homeomorphism.  This means that there is $\epsilon >0$ such that if
\begin{equation}
\sup_{x\in {\Bbb B} } q(f(x),\varphi(x))<\epsilon,
\end{equation} 
then $f$ is not $K$--quasiconformal for any $K<\infty$.
\end{corollary}

The methods developed by Donaldson and Sullivan perhaps allow one to extend the Atiyah--Singer\index{Atiyah--Singer}
index theory of the first order elliptic differential operators to quasiconformal $4$--manifolds (and to other even
dimensions),  and to study the  de Rham cohomology.  Earlier results along these lines had
been developed and studied by Teleman using the Lipschitz structures on topological
$n$--manifolds guaranteed by Sullivan's results.

Recent results of Sullivan and others seem to suggest that the  Seiberg-Witten equations
cannot be used so effectively in the quasiconformal category. Thus there is perhaps a distinction between the
topological, quasiconformal and smooth categories in dimension $4$.

\subsection{The extension problem}\index{quasiconformal extension}

The deformation and approximation theory developed by Sullivan for quasiconformal mappings has other important applications. One of these is the extension or boundary value problem: can a quasiconformal homeomorphism $f:{\Bbb R}^n\to{\Bbb R}^n$ be lifted to a quasiconformal mapping of ${\Bbb R}^{n+1}$.  Actually, since $f(\infty)=\infty$ defines a quasiconformal homeomorphism of the Riemann sphere  ${\Bbb S}^n\approx {\Bbb R}^n\cup \{\infty\}$ the problem is usually formulated as asking if  given a quasiconformal homeomorphism $f:{\Bbb S}^n\to{\Bbb S}^n$,  is there a quasiconformal $F:{\Bbb B}^{n+1}\to{\Bbb B}^{n+1}$ such that $F|{\Bbb S}^{n}=f$ ?  If the answer is ``yes'' we  would also like it to be quantitative.

The answer to this rather elementary question took rather a long time to find.  In one dimension it is a well known theorem of Ahlfors and Beurling concerning the boundary values of quasiconformal mappings of the disk (and quasisymmetric mappings). In dimension three Carleson gave a proof which relied on some combinatorial/piecewise linear topology which is not available in higher dimensions \cite{Carleson}.  Tukia and V\"ais\"al\"a developed and applied Sullivan's ideas to solve this problem,  \cite{TV}.  We remark that it is not at all obvious that this should be the case though.  For instance,  it follows from Milnor's work that there is a diffeomorphism $f:{\Bbb S}^6 \to{\Bbb S}^6$ which cannot be extended to a diffeomorphism of ${\Bbb B}^7$. The solution to the lifting problem shows there is however a quasiconformal extension (since a diffeomorphism of the sphere is quasiconformal).

\begin{theorem}  Let $K\geq 1$.  There is $K^*=K^*(n,K)$ such that if $f:{\Bbb S}^n\to{\Bbb S}^n$ is $K$-quasiconformal, then there exits a $K^*$-quasiconformal $F:{\Bbb B}^{n+1}\to{\Bbb B}^{n+1}$ such that $F|{\Bbb S}^{n}=f$.
\end{theorem}
This theorem is quite nontrivial to establish.  Basically one constructs an obvious extension which will not in general be a homeomorphism but is  ``almost quasiconformal'' at large scales in the hyperbolic metric of the ball ${\Bbb B}^n$.  Such things are called quasi-isometries in the literature.  The approximation results of Sullivan \& Tukia-V\"ais\"al\"a show that such mappings can be approximated in the $C^0$ topology of the hyperbolic metric of the ball by quasiconformal maps.  The technical condition they require is called $\varphi$-solid.  Finally,  any two maps which are a bounded distance apart in the hyperbolic metric agree on the sphere -- an elementary consequence of  hyperbolic geometry  -- and so the constructed quasiconformal approximation is an extension of the given boundary values.

\subsection{Boundary values of quasiconformal mappings}

The converse problem to the problem  discussed above is well known and rather easier:

\begin{theorem}\label{extn}  Let $F:{\Bbb B}^n\to{\Bbb B}^n$ be a $K$-quasiconformal homeomorphism. Then $F$ extends quasiconformally to $F^*:\hat {{\Bbb R}}^n\to\hat {{\Bbb R}}^n$,  $F^*|{\Bbb B}^n=F$.  Further,  if $f=F^*|{\Bbb S}^{n-1}$, $n\geq 3$,  then
\[ f:{\Bbb S}^{n-1}\to{\Bbb S}^{n-1}\]
 is $K$-quasiconformal.  
\end{theorem}
Actually for this theorem we may as well assume $F(0)=0$ and the extension can be effected by reflection;  for $|x|>1$ define
\[ F^*(x) = \frac{F(x/|x|^2)}{|F(x/|x|^2)|^2}. \]
The difficulty now lies in establishing that $F$ has a continuous extension to the boundary,  but this follows from quite direct modulus estimates.  Since $F^{-1}$ also satisfies the same hypothesis,  it also has a continuous extension to ${\Bbb S}^{n-1}$,  so $F$ has a homeomorphic extension and it directly follows that this extension is quasiconformal.

\medskip

This leads one directly to consider the Carath\'eodory problem for the boundary values of quasiconformal mappings.  Recall that Carath\'eodory proved that a conformal map of the unit disk  $\varphi :{\Bbb D}\to\Omega$ extends homeomorphically to the boundary if and only if $\partial\Omega$ is a Jordan curve.  The topological obstructions to such a result in higher dimensions are manifest,  especially considering that the quasiconformal image of a ball could have as boundary a wildly knotted sphere.  However,  it is quite clear the extension result will remain valid if locally the boundary is quasiconformally equivalent to the boundary of the ball (a notion referred to as collaring). However not much beyond this is known.

\medskip

Finally here,  the reader familiar with complex analysis will be well aware of the substantial theory around the structure and properties of bounded analytic functions.  There is an analogous theory for quasiregular mappings and while there are some interesting results there remains some substantial issues to be resolved.  A major question concerns the existence almost everywhere of radial limits,  a well known and useful result for analytic mappings.  ($f$ has radial limits at $\zeta\in \partial{\Bbb B}^n$ if whenever ${\Bbb B}^n\ni x_n\to\zeta$ so that $1-|x_n|\approx |x_n-\zeta|$,  then $f(x_n)$ has a limit).  Thus one might ask the following: given a bounded quasiregular mapping $f:{\Bbb B}^n\to{\Bbb R}^n$ is it true that $f$ has radial limits almost everywhere ?  At this point I believe it is not known even if $f$ has a single radial limit.  These sorts of results are known with additional assumptions such as finite Dirichlet energy,  see eg. \cite{Vu}.

\subsection{Generalised Beltrami systems}
 
Recall that a measurable conformal structure\index{conformal structure!measurable} on a domain $\Omega\subset{\Bbb R}^n$ is a measurable map $G:\Omega \to S(n)$,  the non-positively curved symmetric space of positive definite symmetric $n\times n$ matrices of determinant equal to $1$.   We will always assume that such a map is bounded and such an assertion is equivalent to the assumption that there is a constant $K<\infty$ such that 
\begin{equation}\label{bmcs} \max_{|\zeta|=1} |G(x) \zeta| \leq K\min_{|\zeta|=1} |G(x)\zeta |\end{equation}
for almost every $x\in \Omega$.  The number $K$ plays the role of an ellipticity constant in the associated nonlinear PDE we shall encounter.

We can use these ideas to study Beltrami systems on manifolds.  In what follows we avoid technicalities by simply discussing what happens locally in ${\Bbb R}^n$ - the tangent space to a smooth manifold.

The bounded measurable conformal structure $G$ can be used to define an inner-product on the tangent spaces to $\Omega$ by the rule
\begin{equation}
\langle u,v \rangle_{G} = \langle u,G(x) v \rangle, \hskip30pt u,v\in T\Omega_x.
\end{equation}
Thus (\ref{bmcs}) implies that the unit balls in the metric $\langle , \rangle_{G} $ on the tangent space have uniform eccentricity when viewed in the Euclidean metric $\langle, \rangle$.  For this reason a measurable conformal structure is often referred to as a bounded ellipse field.\\

Suppose now that $\tilde{\Omega}$ is another domain and $H:\tilde{\Omega}\to S(n)$ is a measurable conformal structure defined on $\tilde{\Omega}$.

The {\em generalised Beltrami system} is the PDE
\begin{equation}\label{GHdef}
D^tf(x) \, H(f(x)) \,Df(x) = J(x,f)^{2/n}\, G(x) \hskip20pt \mbox{almost every $x\in\Omega$}
\end{equation}
where a solution $f:\Omega\to\tilde{\Omega}$ is  assumed to be a mapping of Sobolev class $W^{1,n}_{loc}(\Omega)$.  We wish to place an {\em ellipticity condition} on equation (\ref{GHdef}) to link this with the theory of quasiconformal and quasiregular mappings.  This takes the form
\begin{equation}\label{ellipass}
\|d_S(G,I_n)\|_\infty+\|d_S(H,I_n)\|_\infty \leq M <  \infty,
\end{equation}
where $d_S$ is the metric of $S(n)$ and $I_n$ is the $n\times n$ identity matrix.  This metric is discussed in Wolf's book \cite{Wolf}.  The assumption at (\ref{ellipass}) bounds from above and below the ratio of the largest to the smallest eigenvalues of $G$ and $H$ and applying norms shows that (\ref{GHdef}) together with (\ref{ellipass}) gives the existence of a constant $K=K(M,n)$ so that 
\[ |Df(x)|^n \leq K J(x,f),  \]
that is $f$ is quasiconformal.

\medskip

Next,  the following calculation is very informative.  If $u,v\in T\Omega_x$,  then almost everywhere
\begin{eqnarray*}
\langle f_*u, f_*v\rangle_{H} & = & \langle Df(x)u, Df(x) v\rangle_H = \langle Df(x)u, H(f(x)) Df(x) v\rangle_H \\
& = &  \langle  u, D^tf(x)H(f(x)) Df(x) v\rangle  =   \langle  u, J(x,f)^{2/n}\,G(x) v\rangle \\
& = &  J(x,f)^{2/n} \langle  u,   v\rangle_G.  \end{eqnarray*}
This shows that $f$ preserves the inner-product between tangent vectors up to a scalar multiple.   Therefore $Df$ preserves angles between tangent vectors and $f$ the ante between curves (almost everywhere).  Thus $f$ can be viewed as a conformal mapping between the spaces $(\Omega,G)$ and  $(\tilde{\Omega},H)$.

\medskip

It is fair to say the theory of the equation (\ref{GHdef})  is complete and about as good as one could wish for in two-dimensions.  This is because when written in complex notation and with a bit of simplification it reduces to the linear first order equation
\begin{equation}
\frac{\partial f}{\partial \bar z} = \mu(z) \; \frac{\partial f}{\partial  z} + \nu(z) \overline{\frac{\partial f}{\partial  z} } 
\end{equation}  
with the ellipticity bounds $|\mu(z)|+|\nu(z)|\leq k < 1$ for almost all $z$.  The measurable functions $\mu$ and $\nu$ can be explicitly determined from $G$ and $H$.

When $H$ is the  identity,  we have $\nu\equiv 0$ and the usual Beltrami equation,
\begin{equation}
\frac{\partial f}{\partial \bar z} = \mu(z) \; \frac{\partial f}{\partial  z},
\end{equation}  
with $\|\mu\|_\infty = k<1$ which readers familiar with Teichm\"uller theory will no doubt recognise.  A thorough modern account of the theory of Beltrami equations is given in \cite{AIM}.   

\medskip

In higher-dimensions we have already commented above on the various forms of topological rigidity \index{topological rigidity} that occur for solutions to Beltrami systems.  It is basically the following result which assures us we are going to have to deal with discontinuous conformal structures if there is to be a viable theory of branched mappings preserving a conformal structure.  We will call these things ``rational mappings'' later. 

\begin{theorem}  Let $f:\Omega\to f(\Omega)$ be a $W^{1,n}_{loc}(\Omega,f(\Omega))$, $n\geq 3$,  solution to the equation (\ref{GHdef}) where both $G$ and $H$ are continuous.  Then $f$ is a local homeomorphism.
\end{theorem} 

This result can be strengthened in various ways -- continuity is really too strong here, see \cite{IM}. As an easily seen consequence of Theorem \ref{MartioC} we further note that that if the ellipticity constants are close enough to $1$,  that is if $G$ and $H$ are sufficiently and uniformly close to the identity,  then $f$ is also a local homeomorphism.

\section{Nevanlinna Theory}

The classical theorem of Picard  of 1879 initiated the value distribution theory of holomorphic functions in the complex plane.  It simply states
that an entire function which omits two values is constant.  Nevanlinna theory  
is a far reaching extension of Picard's theorem and concerns the distribution of the values of an
entire function.  It was developed
around 1925.  Ahlfors subsequently brought many new geometric ideas,  including the use of quasiconformal mappings,  into Nevalinna theory.    Given a meromorphic function $f:{\Bbb C}\to {\hat{{\Bbb C}}}$ we
define for any Borel set $\Omega$ and any $y\in {\Bbb C}$ the counting
function
\begin{equation}
n(\Omega,y) = \#\{ f^{-1}(y) \cap \Omega \}
\end{equation}
where the number of points is counted according to multiplicity.  The function $A(r)$ is defined to be the average
of $n(r,y)=n({\Bbb B}(r),y)$ with respect to the spherical measure on ${\hat{{\Bbb C}}}$.  An important result in the area is Ahlfors'
theorem concerning the so--called {\em defect relation}. Given a nonconstant meromorphic function there is a set
$E\subset [1,\infty)$ of finite logarithmic measure
\begin{equation}
\int_{E} \frac{dr}{r} < \infty
\end{equation}
such that for $a_1,a_2,\ldots a_q$ distinct points in ${\hat{{\Bbb C}}}$,
\begin{equation}
\limsup_{E \not\ni r \to \infty} \sum_{j=1}^{q} \delta(a_j,r) \leq 2
\end{equation}
where 
\begin{equation}
\delta(a_j,r) = \max \{1-\frac{n(r,a_j)}{A(r)},0\}
\end{equation}
is called the defect of $a_j$.  Roughly,  for $r$  off a thin set $E$
the function $f$ covers each point $a_j$ the correct``average number'' of times on the ball of radius $r$.  Picard's
\index{Picard} theorem is a direct consequence if we put in the three omitted values for the numbers $a_j$.  There is of
course much more to value distribution theory in the plane than this.  However,  here we would like to mention the
$n$--dimensional analogues of these results for quasiregular mappings.

In a series of brilliant papers  presented in his monograph
\cite{R}, S. Rickman developed the value distribution theory of
quasiregular mappings using geometric methods for the most part,  and
in particular extremal length.  The most striking result so far
obtained is the sharp form of the defect relation. (The definitions of
counting functions and so forth in higher dimensions are the obvious generalisations).\index{Rickman's theorem}

\begin{theorem}[Rickman's theorem]  There is a constant $C(n,K)$ such that if $f:{\Bbb R}^n \to \hat {{\Bbb R}}^n$ is a nonconstant $K$--quasiregular mapping
there is a set $E$ of finite logarithmic measure such that
\begin{equation}
\limsup_{E \not\ni r \to \infty} \sum_{j=1}^{q} \delta(a_j,r) \leq 2
\end{equation}
whenever $a_1,a_2,\ldots a_q$ are distinct points of $\hat {{\Bbb R}}^n$. 
\end{theorem}
As a consequence we obtain the following version of the Picard Theorem
\begin{theorem}
For each $K\geq 1$ there is an integer $q=q(n,K)$ such that every $K$--quasiregular map
$f:{\Bbb R}^n \rightarrow \hat {{\Bbb R}}^n \backslash \{ a_1,a_2,\ldots,a_q \}$,  where $a_j$ are distinct, is
constant. 
\end{theorem}
It was thought for a while that the number $q(n,K)=2$. However, Rickman gave an example to show that
this is not the case,  at least when $n=3$.
\begin{theorem}
For every positive integer $p$ there exists a nonconstant $K$--quasiregular mapping  $f:{\Bbb R}^3
\rightarrow {\Bbb R}^3$ omitting $p$ points. 
\end{theorem}
Such examples are highly non-trivial to construct.  However the theory needs more such examples and a better understanding of what is going on here.  And although Rickman's result is expected to hold in dimensions $n\geq 4$,  this has not yet been confirmed.  Lewis \cite{Lewis},  following joint work with Eremenko \cite{EL}, gave an analytic proof of the quasiregular version of the Picard
Theorem using non-linear potential theory and in particular Harnack's inequality for
$\cal A$--harmonic functions. This proof is refined in \cite[Chapter17]{IM}.

\medskip

Next we recall Rickman's version of Montel's theorem.  This result also has been generalised further to consider quasiregular mappings into manifolds with a suitable number of ends (which play the role of omitted points).\index{Montel's theorem}

\begin{theorem}[Montel's theorem]\label{Montel} For every $K\geq 1$ and $n\geq 2$ there is an integer $q_{n,K}$ with the following property.  Let $\epsilon>0$.  If  ${\cal F}$ is a family of $K$--quasiregular mappings $f:\Omega\subset\hat {{\Bbb R}}^n\to\hat {{\Bbb R}}^n$ such that each $f\in {\cal F}$ omits $q$ values $a_{1}^{f},a_{2}^{f},\ldots,a_{q_{n,K}}^{f}$ for which the spherical distances
\[ \sigma(a_{i}^{f},a_{j}^{f}) > \epsilon, \hskip20pt i\neq j,\]
then ${\cal F}$ is a normal family.
\end{theorem}
Note here that $q_{n,K}$ does not depend on $\epsilon$.  Then an elementary topological argument gives a sharper version in the quasiconformal case.

\begin{theorem} Let $\epsilon>0$.  If  ${\cal F}$ is a family of $K$--quasiconformal mappings $f:\Omega\subset\hat {{\Bbb R}}^n\to\hat {{\Bbb R}}^n$ such that each $f\in {\cal F}$ omits two values $a_{1}^{f}$ and $a_{2}^{f}$ for which the spherical distance $\sigma(a_{1}^{f},a_{2}^{f}) > \epsilon$,  then ${\cal F}$ is a normal family.
\end{theorem}

Another useful normal families criterion is through Zalcman's lemma in higher dimensions due to Miniowitz \cite{Miniowitz}.   We say a family of mappings ${\cal F}$ is normal at a point $x_0$ if there exists an open neighbourhood $U$ of $x_0$ on which the family ${\cal F}|U=\{f|U:f\in {\cal F}\}$ is normal.

\begin{theorem}[Zalcman's Lemma] \label{zalcman}\index{Zalcman lemma} Let $K\geq 1$ and ${\cal F}$ a family of $K$-quasiregular mappings $f:{\Bbb B}^n\to\hat {{\Bbb R}}^n$.  Then ${\cal F}$ is not normal at $x_0\in {\Bbb B}^n$ if and only if there is a sequence of positive numbers $r_j\searrow 0$,  a sequence of points $x_j\to x_0$ and a sequence of mappings $\{f_j\}_{j=1}^{\infty}\subset{\cal F}$ such that if we define
\begin{equation}
\varphi_j(x)  = f_j(x_j+r_j x),
\end{equation}
then $\varphi_j$ converges uniformly on compact subsets of ${\Bbb R}^n$ to a non-constant $K$--quasiregular mapping $\varphi:{\Bbb R}^n\to\hat {{\Bbb R}}^n$.
\end{theorem}
The term quasimeromorphic is sometimes used for quasiregular mappings which assume the value $\infty$ continuously in the spherical metric.   Again,  these theorems and their near relatives remain true in much more general settings.  These have been worked out by Rickman,  his students and others,  see for instance \cite{HR} and the references therein.

\section{Non-linear Potential Theory}

The nonlinear potential theory as it pertains to the theory of higher-dimensional quasiconformal and quasiregular mappings is quite comprehensively covered in the book of Heinonen, Kilpal\"ainen and Martio \cite{HKM}.  One of the central results is the fact that quasiconformal and quasiregular mappings are morphisms for the class of $p$-harmonic functions.  Typically the modern theory deals with structurally nice measures $d\mu$ defined on  ${\Bbb R}^n$ and   giving weighted Sobolev space.  These measures,  called {\em admissible},  are almost always of the form $d\mu=\omega(x)\; dx$, so
\[ \mu(E) = \iint_E \omega(x) \; dx \]
 and $\omega$ is an admissible weight satisfying the four conditions
\begin{enumerate}
\item {\bf Doubling:} There is a constant $c_1$ such that
\[ \mu(2B) \leq c_1\mu(B),  \hskip20pt B={\Bbb B}^n(x,r).\]
\item {\bf Testing:}  If $\Omega$ is open and $\varphi_i\in C^{\infty}(\Omega)$ with 
\[ \iint_\Omega|\varphi_i|^p \; d\mu \to 0, \hskip10pt {\rm and}\hskip10pt \iint_\Omega|\varphi_i -  v|^p \; d\mu \to 0\]
for vector valued $v$ in $L^{p}_{\mu}(\Omega)$,  then $v=0$.
\item {\bf Sobolev Embedding:} There are constants $\alpha>1 $ and $c_2$ such that for all balls $B={\Bbb B}(x_0,r)\subset{\Bbb R}^n$ and all $\varphi\in C^{\infty}_{0}(B)$ we have
\[ \Big( \frac{1}{|\mu(B)|} \iint_B |\varphi|^{\alpha p} \; d\mu \Big)^{1/(\alpha p)} \leq c_2 r \Big(\frac{1}{|\mu(B)|} \iint_B |\varphi|^{ p} \Big)^{1/p} \; d\mu. \]
\item {\bf Poincar\'e Inequality:} There is a constant $c_3$ such that for all balls $B={\Bbb B}(x_0,r)\subset{\Bbb R}^n$ and all bounded $\varphi\in C^{\infty}(B)$ we have
\[   \iint_B |\varphi-\varphi_B|^{ p} \; d\mu \leq c_3 r^p \;  \iint_B |\nabla \varphi|^{ p} \; d\mu \]
\noindent where $\varphi_B$ is the $\mu$-average of $\varphi$,  $\varphi_B=|\mu(B)|^{-1} \iint_B \varphi \; d\mu $.
\end{enumerate}

\medskip

The constants $c_1,c_2$ and $c_3$ do not matter so much,  but the constant $\alpha$ plays an important role in regularity.  From these one directly gets the weighted Poincar\'e inequality\index{Poincar\'e inequality}
\begin{theorem}[Poincar\'e inequality]  If $\Omega$ is a bounded domain,  then for all $\varphi\in C^{\infty}_{0}(\Omega)$ we have
\[ \iint_\Omega |\varphi|^p \;d\mu \leq c_{2}^{p} \;{\rm diam}^p(\Omega) \iint_\Omega |\nabla\varphi|^p \; d\mu \]
\end{theorem}
The first direct connection with the theory of higher-dimensional quasiconformal mappings is the following.
\begin{theorem}  If $f:{\Bbb R}^n\to{\Bbb R}^n$ is quasiconformal and $J(x,f)$ denotes its Jacobian determinant,  then
\[ \omega(x) = J(x,f)^{1-p/n} \]
is an admissible weight whenever $1<p<n$.
\end{theorem}
There is now a substantial literature on admissible weights and their role in generalising many of the basic results of analysis to more general settings (including metric spaces).

\subsection{${\cal A}$-harmonic functions}\index{${\cal A}$-harmonic functions}

One of the key tasks of nonlinear potential theory is to develop techniques to study the quasilinear elliptic equation  
\begin{equation} \label{genA} - {\rm div} \; {\cal A}(x,\nabla u) = 0, \end{equation}
generalising the Laplace equation ${\cal A}(x,\zeta)=\zeta$.  To get any viable theory one needs structural conditions on ${\cal A}$ and these are usually of the form
\[ {\cal A}(x,\zeta)\cdot\zeta \approx \omega(x) |\zeta|^p \]
for an admissible weight $\omega$,  together with various (less important) technical assumptions.  The weighted $p$-Laplace equation is
\begin{equation}
- {\rm div} (\omega(x)|\nabla u|^{p-2} \nabla u )= 0,
\end{equation}
the solutions of which are local minimisers of the weighted energy integral
\[ \iint |\nabla u|^p \; \omega(x) \; dx. \]
In a similar manner,  the general ${\cal A}$-harmonic equation is connected with the local extrema,  satisfying the Euler-Lagrange equations,  of variational integrals of the form
\[ \iint F(x,\nabla u)\; dx.  \]

\medskip

The precise assumptions on ${\cal A}$ take the following form.  Suppose that $0<\alpha\leq \beta <\infty$ and suppose ${\cal A}:{\Bbb R}^n\times{\Bbb R}^n\to{\Bbb R}^n$ is a mapping such that $x\mapsto {\cal A}(x,\zeta)$ is measurable for all $\zeta\in {\Bbb R}^n$, and $\zeta\mapsto {\cal A}(x,\zeta)$ is continuous for almost every $x\in {\Bbb R}^n$ and finally we require that for all $\zeta\in{\Bbb R}^n$ and almost all $x\in {\Bbb R}^n$ we have the estimates:
\begin{itemize}
\item ${\cal A}(x,\zeta)\cdot\zeta \geq \alpha \omega(x) |\zeta|^p$,
\item $|{\cal A}(x,\zeta)|\leq \beta \omega(x)|\zeta|^{p-1}$,
\item $({\cal A}(x,\zeta_1)-{\cal A}(x,\zeta_2))\cdot(\zeta_1-\zeta_2) > 0,  \hskip15pt  \zeta_1\neq\zeta_2$,
\item ${\cal A}(x,\lambda \zeta)=\lambda |\lambda|^{p-2} {\cal A}(x,\zeta), \hskip15pt  \lambda\in{\Bbb R}\setminus\{0\}$.
\end{itemize}

These conditions describe the allowable degenerate behaviour of the equation (\ref{genA}).  A function $u\in W^{1,p}_{loc}(\Omega)$ is a {\em weak solution} to the equation (\ref{genA}) in $\Omega$ if for every $\varphi\in C^{\infty}_{0}(\Omega)$ we have
\begin{equation}
\iint_\Omega {\cal A}(x,\nabla u(x))\cdot \nabla \varphi(x) \; dx = 0,
\end{equation}
and a {\em super-solution} if the left-hand side is non-negative.   A real-valued function $h:\Omega\to {\Bbb R}$ is called ${\cal A}$-harmonic if it is a continuous weak solution to (\ref{genA}).  These functions are the main object of study in the theory in as much as the harmonic functions are for classical potential theory.

We then have the following theorem telling us the Dirichlet problem has a solution for Sobolev boundary values.

\begin{theorem} Let $\Omega\subset {\Bbb R}^n$ be a bounded domain and $\eta \in W^{1,p}_{\mu}(\Omega)$.  Then there is a unique solution $u\in W^{1,p}_{\mu}(\Omega)$ of the equation (\ref{genA}) with $u-\eta\in W^{1,p}_{0,\mu}(\Omega)$.
\end{theorem}
The notation should be self explanatory,  but the last space here consists of those functions whose $p^{th}$ power of their absolute value is integrable with respect to $d\mu$ and which vanish at the boundary in the Sobolev sense.

\medskip  

The theory then develops by studying regularity (giving continuity) and compactness properties of solutions (locally uniformly bounded families are equicontinuous), firmly establishing the connection with the variational formulation and then the Harnack principle and maximum principle.  For instance
\begin{theorem}[Strong maximum principle] A non-constant ${\cal A}$-harmonic function defined in a domain $\Omega$ cannot achieve its maximum or minimum value.
\end{theorem}
Of course once connected with quasiconformal mappings the maximum principle will imply that quasiregular mappings are open.
\medskip

The capacity and ${\cal A}$-harmonic measure theory are deep and interesting and have important consequences,  but it would lead us too far astray to develop them here.

\subsection{Connections to quasiconformal mappings}

If $f=(f^1,f^2,\ldots,f^n)$ is a solution to the Beltrami system (\ref{BE}),  then  $u=f^i$
satisfies the following equation of  elliptic type:
\begin{equation}\label{Aharm}
{\rm div}(\langle G^{-1}\nabla u, \nabla u \rangle ^{(n-2)/2} G^{-1} \nabla u ) =0,
\end{equation}
where $G^{-1}=G^{-1}(x)$ is the inverse of the distortion tensor.  This is easily seen by unwinding the tautological Beltrami equation that $f$ satisfies.

Therefore  $u$  is an $\cal A$-harmonic
function.  With the choice $G={\rm Id}_{n\times n}$ (the Cauchy-Riemann system and conformal mappings) we see that $u$ is $n$-harmonic.  These ideas lead to an alternative proof of the Liouville theorem which was developed by Reshetnyak.

Actually,  if  $f$  is a solution of the
Beltrami system,  then  remarkably $u=\log |f|$  satisfies an $\cal A$--harmonic equation as well.  

\begin{theorem}  Let $f:\Omega\to{\Bbb R}^n$ be a non constant quasiregular mapping and let $b\in {\Bbb R}^n$.  Then the function
\[ u(x) = \log |f(x)-b| \]
is ${\cal A}$-harmonic in the open set $\Omega\setminus f^{-1}(b)$.  Here ${\cal A}$ satisfies the structure equations with $p=n$,  $\omega(x)=1$,  $\alpha=1/K$ and $\beta=K$.
\end{theorem}

Now the  theory of  $\cal A$-harmonic functions implies directly that quasiregular
mappings are open and discrete (we discussed open mapping property above).  Roughly,  discreteness follows from the omitted discussion of ${\cal A}$-harmonic measure.  Here we need
the fact that the polar sets ($\{ x: \log |f(x)-b| = -\infty \}$) have conformal capacity zero.  This implies they have Hausdorff dimension zero and hence are totally disconnected.  A topological degree argument then completes the proof of Reshetnyak's theorem -  discussed earlier  as Theorem \ref{RST}.

\subsection{Removable singularities} \index{removable singularity}

There is a classical theorem of Painlev\'e concerning removable sets for analytic functions.  It
states that if $\Omega \subset {\Bbb C}$ is a planar domain,  $E\subset \Omega$ is a closed subset and
$f:\Omega \backslash E \rightarrow {\Bbb C}$ is a bounded analytic function,  then  $f$ has an analytic
extension to $\Omega$.  This result has also found generalisation in higher dimensions and we give a brief account of that here.  Details are to be found in \cite{IM}

\medskip

A closed set $E\subset {\Bbb R}^n$ is {\em removable under bounded $K$--quasiregular
mappings} if for every open set $\Omega \subset {\Bbb R}^n$ any bounded $K$--quasiregular mapping 
$f:\Omega \backslash E \rightarrow {\Bbb R}^n$ extends to a $K$--quasiregular mapping of $\Omega$.  We
stress here that $f$  need not even be locally injective,  nor even of bounded topological degree. 
 
\begin{theorem}
There is $\epsilon=\epsilon(K)>0$ such that closed sets of Hausdorff dimension $\epsilon$ are removable under bounded $K$--quasiregular mappings.
\end{theorem}
In particular,  sets of Hausdorff dimension $0$ are always removable for bounded quasiregular mappings.  In light of 
conjectures regarding the
$p$--norms of the Hilbert transform on forms and the relationship between $s$--capacity and Hausdorff dimension  we  
formulate the following conjecture regarding the optimal result,  see \cite{IM}.
\begin{conjecture}
Sets of Hausdorff $d$--measure zero,  $d=n/(K+1) \leq n/2$, are removable under bounded
$K$--quasiregular mappings. 
\end{conjecture}
In $2$--dimensions Astala \cite{Astala} has   verified this conjecture for all $d<2/(K+1)$ and the borderline cases are well in hand.
In response to these questions Rickman \cite{Rickman3} has constructed examples to
show that the results are,  in some sense, best possible. 

\begin{theorem} There are Cantor sets $E$ of arbitrarily small
Hausdorff dimension and bounded quasiregular mappings ${\Bbb R}^3 \backslash E \rightarrow {\Bbb R}^3$.  For
such mappings  $E$  is necessarily non-removable. \end{theorem}

Here we must have $K\to\infty$ as the Hausdorff dimension tends to zero.

\section{Quasiregular   dynamics in higher dimensions}

There are a number of recent developments in the theory of higher-dimensional quasiconformal and quasiregular mappings which link these areas to questions of dynamics and different types of rigidity phenomena in dimension $n\geq 3$ than those we have seen earlier.  A related survey can be found in \cite{Mart2}.  A self-mapping of an $n$-manifold is {\em rational} or {\em uniformly quasiregular}  if it preserves some bounded measurable conformal structure (see below at (\ref{GG})).  In what follows we will assume that the manifold in question is Riemannian for simplicity,  but Sullivan's theorem shows that the existence of a bounded measurable Riemannian structure is a purely topological notion, at least when $n\neq 4$.  The bounded measurable structure which is preserved will typically {\em not} be the underlying Riemann structure -- indeed most often this structure will necessarily be  discontinuous.

There is a close analogy between the dynamics of rational maps of closed manifolds and the classical Fatou-Julia theory of iteration of rational mappings of ${\hat{{\Bbb C}}}$.  The theory is particularly interesting on the Riemann $n$-sphere $\hat {{\Bbb R}}^n$ where many classical results find their analogue,  some of which we will discuss below.  In higher dimensions other interesting aspects come into play.  An analytic mapping of a closed surface is a  homeomorphism unless the surface is  ${\Bbb S}^2$ or the $2$-torus - where of course the map is covered by multiplication on the complex plane.    Thus there cannot be any interesting Fatou/Julia type theories on surfaces other than ${\Bbb S}^2$ where it is very highly developed.  We cannot expect this situation to persist exactly in higher dimensions,  but informed by the two-dimensional case we might expect an interaction between the curvature of the manifold and the existence or otherwise of rational mappings. Thus we introduce the Lichnerowicz problem of classifying those manifolds admitting rational endomorphisms.   Once we have examples of nontrivial rational maps of such spaces (for instance the $n$-spheres) we can ask to what extent the classical Fatou/Julia theory remains true.  What is the structure of the Julia set and can we classify dynamics on the Fatou set.

\medskip

Recall that measurable conformal structure on $\hat {{\Bbb R}}^n$ is a measurable map $G:\hat {{\Bbb R}}^n \to S(n)$,  the non-positively curved symmetric space of positive definite symmetric $n\times n$ matrices of determinant equal to $1$.  A $W^{1,n}(\hat {{\Bbb R}}^n)$ map 
\[ f:\hat {{\Bbb R}}^n\to\hat {{\Bbb R}}^n\]
 will preserve this conformal structure if it satisfies the generalised Beltrami system
\begin{equation}\label{GG}
D^tf(x) \, G(f(x)) \,Df(x) = J(x,f)^{2/n}\, G(x) \hskip20pt \mbox{almost every $x\in\hat {{\Bbb R}}^n$}.
\end{equation}
With an ellipticity assumption on $G$ as before,  there is a $K<\infty$ such that any $W^{1,n}(\hat {{\Bbb R}}^n)$ solution to (\ref{GG}) is $K$-quasiregular.  The composition of quasiregular maps is certainly quasiregular,  and so the set of solutions to the equation (\ref{GG}) forms a semigroup under composition.  Let us denote the semigroup of solutions to (\ref{GG}) by ${\rm Rat}(G)$,  the maps rational with respect to the measurable conformal structure $G$. 

  \subsection{Existence of equivariant measurable conformal structures}  
  
The nonpositive curvature of the space $S(n)$ allows one to make various averaging constructions, first noted by Sullivan in two-dimensions for group actions and develop by Tukia in higher-dimensions.  Slightly refining these arguments gives us the following theorem,  see \cite{IM}[Chapter 21].
  
  \begin{theorem}[Semigroup]\label{semigroup}   Let ${\cal F}$ be an abelian semigroup of quasiregular mappings of a manifold $M^n$ such that each $f\in{\cal F}$ is $K$-quasiregular.  Then there is a measurable conformal structure $G_{\cal F}:M^n\to S(n)$ and for each $f\in{\cal F}$, 
  \begin{equation} 
D^tf(x) \, G_{\cal F}(f(x)) \,Df(x) = J(x,f)^{2/n}\, G_{\cal F}(x) \hskip20pt \mbox{almost every $x$}.
\end{equation}
  \end{theorem}
The condition that the semigroup be abelian is too strong,  but some condition is necessary even in two-dimensions, see \cite{Hink}.  In \cite{IM} a ``left-right'' coset condition is given which is automatically true in the abelian case.  In the case of quasiconformal groups,  the result is more general,  but because of  Lelong's Theorem \ref{lich} (described below) the result is really only of interest in ${\Bbb B}^n$,  ${\Bbb R}^n$ and $\hat {{\Bbb R}}^n$.

  \begin{theorem}[Group] \label{qcgroup}  Let ${\Gamma} $ be a group of quasiconformal self-homeomorphisms of a domain $\Omega \subset \hat {{\Bbb R}}^n$ such that each $g\in{\Gamma} $ is $K$-quasiconformal.  Then there is a measurable conformal structure $G_{{\Gamma} }:\Omega\to S(n)$ and for each $g\in{\Gamma} $, 
  \begin{equation} 
D^tg(x) \, G_{{\Gamma} }(g(x)) \,Dg(x) = J(x,g)^{2/n}\, G_{{\Gamma} }(x) \hskip20pt \mbox{almost every $x\in \Omega$}.
\end{equation}
  \end{theorem}

\subsection{Fatou and Julia sets}

The semigroup of solutions to (\ref{GG}) has the property that the composition of its elements cannot increase the distortion beyond a uniform bound.  This is surprising inasmuch as  if a solution is branched,  then its iterates $f$,  $f\circ f$, \ldots, $f\circ f \circ \cdots \circ f$ have ever increasing degree.  We will see soon that there can even be fixed points of a rational mapping which are also branch points.  It would seem that there should be little or no distortion at a fixed point in order for it not to grow under iteration, yet we know from higher-dimensional topological rigidity that there must be some distortion at branch points.  It is very interesting to see how these observations reconcile.  First we note the following obvious fact.

\begin{theorem}
${\rm Rat}(G)$,  the space of quasiregular solutions to (\ref{GG}),  is closed under composition and there is a $K<\infty$ such that each $f\in {\rm Rat}(G)$ is $K$-quasiregular.
\end{theorem}
Because of Rickman's version of Montel's Theorem (Theorem \ref{Montel})  there is a reasonably complete Fatou-Julia theory associated with the iteration of rational mappings. But first we need to state that there are examples.

\begin{theorem}  For each $n\geq 2$ there is a $K$-quasiregular mapping $f:\hat {{\Bbb R}}^n\to \hat {{\Bbb R}}^n$ with the property that all the iterates $f^{\circ 2}=f\circ $, \ldots, $f^{\circ n}=f\circ f\circ \cdots \circ f^n$ are also $K$-quasiregular.  Thus the family $\{f^{\circ n}:n\geq 1\}$ is a quasiregular semigroup.
\end{theorem}

As a corollary from the existence invariant conformal structures,  Theorem \ref{semigroup},   we have the following.

\begin{corollary}  For each $n\geq 2$ there is a bounded  measurable conformal structure $G$ defined on ${\Bbb S}^n$ which admits non-injective rational mappings.  In these cases  ${\rm Rat}(G)$ is infinite,  contains mappings of arbitrary high degree  and is not precompact.\end{corollary}

\medskip

The study of these quasiregular rational mappings for ${\Bbb S}^{n}\approx \hat {{\Bbb R}}^n$ was initiated in joint work with Iwaniec \cite{IMrat},  but have since been developed by V. Mayer, K. Peltonen and others,  see \cite{MM, MMP, Mayer1, Mayer2}. There are strong restrictions on the geometry and topology of closed manifolds admitting nontrivial rational mappings, for instance  they cannot be negatively curved.  

\medskip

The Fatou set ${\cal F}(f)$\index{Fatou set} of a rational mapping $f$ is the open set where the iterates form a normal family (that is have locally uniformly convergent subsequences).
  The Julia set ${\cal J}(f)$\index{Julia set} is the complement of the Fatou set 
  \[ {\cal J}=\hat {{\Bbb R}}^n\setminus {\cal F}. \]  If the degree of $f\geq 2$,  the only interesting case for us,  then the Julia set is nonempty, closed and  a completely invariant set,
  \[ f^{-1}({\cal J}) = {\cal J}. \] 
  Known examples of Julia sets include Cantors sets,  $\hat {{\Bbb R}}^n$ itself (the Latt\`es type examples),  codimension one spheres and somewhat more complicated sets which separate ${\Bbb R}^3$ into infinitely many components.  There are very many interesting and unanswered questions about what sets could be Julia sets.  For instance it is known in three-dimensions,  that only very simple knots (torus knots in ${\Bbb S}^3$) can possibly be Julia sets.

  \subsection{Dynamics of rational mappings}

 We first consider the classification of fixed points.
 In \cite{HMM} it was shown that uniformly quasiregular mappings are locally Lipschitz near a fixed point $x_0$ which is not a branch point. This is then used to show that the family ${\cal F}=\{f_\lambda\colon \lambda>1\}$ is a normal family, where $f_\lambda(z)=\lambda f(z/\lambda)$. This is relatively straightforward as
\[ f_\lambda \circ f_\lambda = (f\circ f)_\lambda \]
and so $f^{\circ 2}_{\lambda}$ is Lipschitz and linearizes again.  Moreover, this more or less implies that all limits of convergent subsequences of ${\cal F}$ are uniformly quasiconformal mappings -- that is the cyclic group $\langle g \rangle=\{g^n:n\in {\Bbb Z}\}$ for $g\in {\cal F}$ is a uniformly quasiconformal group as discussed in the next section. The set of all such limit mappings is called the generalized derivative of $f$ at $x_0$. Uniformly quasiconformal mappings have been classified as either loxodromic, elliptic or parabolic. It follows that the elements of the generalized derivative are either all constant, all elliptic, or all loxodromic, and this allows for a classification of the fixed points of a uniformly quasiregular mapping as attracting/repelling (generalized derivative loxodromic), neutral (generalized derivative elliptic), or super-attracting (fixed point is a branch point).  All these types of fixed points can occur.  Perhaps the most surprising is Mayer's construction of the Latt\`es example.  This example is derived from a functional equation in a similar manner as the classical example,  building on a classification of Martio and Srebro for automorphic quasimeromorphic mappings \cite{MSreb}.  The Julia set in $\hat {{\Bbb R}}^n$ is the unit sphere and both the origin and $\infty$ are super attracting fixed points -  the branch set being a family of lines passing through the origin.  We record this in the following theorem:

\begin{theorem} The Fatou set of a rational quasiregular self-mapping $f$ of $\hat {{\Bbb R}}^n$ have precisely the same types of stable components $U$ as rational functions of $\hat {{\Bbb R}}^2$. They are either
\begin{enumerate}
\item {\bf attracting};  there is $x_0\in U$  with $f(x_0)=x_0$ and as $n\to\infty$,
\begin{equation}\label{goestox0} f^{\circ n}(x) \to x_0 \hskip20pt \mbox{locally uniformly in $U$}\end{equation}
\item {\bf super-attracting};  there is $x_0\in U$  with $f(x_0)=x_0$ and $x_0\in B_f$.  Necessarily then (\ref{goestox0}) holds.   
\item {\bf parabolic};  there is $x_0\in \partial U$  with $f(x_0)=x_0$ and (\ref{goestox0}) holds.
\item {\bf Siegel};  $f:U\to U$ is quasiconformal and $\overline{\{ f^{\circ n} :n\in {\Bbb Z}\}}$ is a compact Lie group.
\end{enumerate}
Further,  there are examples of  types $(1)$, $(2)$ and $(3)$. \end{theorem}

For attracting and repelling fixed points we know the following:

\begin{theorem} A rational quasiregular mapping is locally  quasiconformally conjugate to  the map $x\mapsto 2x$ near a repelling fixed point and is quasiconformally conjugate to $x\mapsto \frac 12 x$ near an attracting fixed point. 
\end{theorem}

Here,  by quasiconformal conjugacy we mean there is a quasiconformal homeomorphism $\varphi$ defined in a neighbourhood $V$ of the origin so that
\[ \varphi\circ f \circ \varphi^{-1} (x) = 2 x \]
Roughly,  after a quasiconformal change of coordinates we get standard dynamics.  There is no standard conformal model for a super-attractive fixed point - a consequence of the Liouville theorem.   This sort of dynamics is quite novel.

There are examples with parabolic dynamics. It can be shown that a map with a parabolic fixed point can be constructed in such a way that it does not admit a quasiconformal linearization in its attracting parabolic petal (unlike the rational case) due to the existence of wild translation arcs.  This builds on Mayer's work in \cite{Mayer3}.

\medskip

\noindent{\bf Question:} An interesting problem is to decide whether or not it is possible to have a ``Siegel disk'' of type $(4)$ described above for a non-injective rational quasiregular mapping.  Such a domain would presumably a ball or solid torus with irrational rotational dynamics.  That the map is injective on a Siegel domain,  and that it would generate a compact Lie group as described does follow.

\medskip

 A classical result is the density of repellors, that is repelling fixed points,  in the Julia set. This is not known in complete generality yet for higher-dimensional rational  mappings,  but we do know the following.
 
 \begin{theorem}  The set of repelling and neutral fixed points is dense in the Julia set.
 \end{theorem} 
 
 In certain cases,  when assumptions on the topological structure of the Julia set,  for instance separating, we do know repellors are dense.

\subsection{Sto\"{i}low Factorisation.}

The following factorisation theorem shows that in fact rational quasiregular mappings are quite common. This is a  variant of the well-known theorem of Sto\"{\i}low's, see \cite{MP2}.

\begin{theorem} \label{theo 1}
Suppose $g:\hat {{\Bbb R}}^n\rightarrow\hat {{\Bbb R}}^n$ is a non-constant quasiregular mapping, $n\ge2$. 
Then there exists  a rational 
quasiregular mapping $f$ whose Julia set is a Cantor set,  and a  quasiconformal  mapping
$h:\hat {{\Bbb R}}^n\to \hat {{\Bbb R}}^n$ such that  $g=f \circ h$.
\end{theorem}

Classically the factorization (for quasiregular maps of ${\hat{{\Bbb C}}}={\mathbb S}^2$ ) is unique up to M\"obius transformation.  If $\varphi \circ f = \psi \circ g$,  then there is a M\"obius transformation  $\Phi$ so that $ \varphi \circ  \Phi= \psi$. Clearly this statement cannot hold in higher dimensions if $\varphi$ and $\psi$ are merely assumed rational.  However if we fix the invariant conformal structure,  then we can make uniqueness statements up to a finite-dimensional Lie group.  Notice the following easy implication showing that no distinction can be made between the branch sets of rational quasiregular mappings an completely general quasiregular mappings, \cite{Mart4}.   

\begin{theorem}  Let $X=B_f$ be the branch set of a quasiregular mapping $f:\hat {{\Bbb R}}^3\to\hat {{\Bbb R}}^3$.  Then $X$ is the branch set of a rational quasiregular mapping.
\end{theorem}

This is a little surprising given how complicated these branch sets can be.  We recall from \cite{HeinRick} that the branch set of a quasiregular mapping could be as wild as Antoine's necklace,  a Cantor set in ${\Bbb S}^3$ whose complement is not simply connected.

\subsection{Smooth rational quasiregular mappings.}
  
The technique used in construction of the factorisation is sufficiently robust that if  the map $f$ is smooth of class $C^k(\hat {{\Bbb R}}^n)$,  then  the quasiconformal homeomorphism $h$,  and consequently the rational mapping $\varphi$,  can be chosen to be smooth of the same class.  Typically one does not expect branched (not locally injective) quasiregular mappings to be smooth,  however there are examples of M. Bonk and J. Heinonen \cite{BH} of a quasiregular map $f:{\mathbb S}^3 \to {\mathbb S}^3$ which is $C^{3-\epsilon}({\mathbb S}^3)$ for every $\epsilon >0$.  Kaufman,  Tyson and   Wu extended these results  to higher dimensions, \cite{KTW}.  The following theorem (which was certainly known to Bonk and Heinonen) is a consequence.

\begin{theorem} There are smooth rational quasiregular mappings of $\hat {{\Bbb R}}^n$ with nonempty branch set, $B_f\neq \emptyset$.  Indeed,
\begin{itemize} 
\item For each $\epsilon>0$,  there is a $C^{3-\epsilon}({\mathbb S}^3)$ rational 
quasiregular mapping $\varphi$ whose Julia set is a Cantor set.

\item For each $\epsilon>0$,  there is a $C^{2-\epsilon}({\mathbb S}^4)$ rational
quasiregular mapping $\varphi$ whose Julia set is a Cantor set.

\item For each $n\geq 5$ there is an  $\epsilon=\epsilon(n)>0$ and a $C^{1+\epsilon}({\mathbb S}^n)$ rational quasiregular mapping $\varphi$ whose Julia set is a Cantor set.
\end{itemize}
\end{theorem}
 Note that although these maps are smooth,  any invariant conformal structure must be quite irregular (at least discontinuous) near the branch set.

 Indeed it is known that the rational mappings described here are structurally stable (or generic);  there is a single attracting fixed point, no relations between critical points and the Julia set is ambiently quasiconformally equivalent to the middle thirds Cantor set.

\subsection{The Lichnerowicz problem: Rational Maps of Manifolds}
 
A natural question now is to ask what sort of manifolds support rational quasiregular endomorphisms. As we have noted,  in two-dimensions it is an easy application of the signature formula for branched coverings to see that only the sphere and torus admit branched self-maps.  In higher dimensions the question is more complicated - though a complete answer was given by  Kangaslampi in three dimensions \cite{Ka}.

\medskip

Here we shall consider such mappings acting on {\em closed} manifolds $M$ of dimension at least two and our problem is   to determine what kind of rational mappings can act on a given manifold. Recall that a rational map $f : M \to M$ is surjective since the continuity and openness of a quasiregular map implies that the image  $fM$ is both compact and open; hence $fM = M$.

\medskip

The first part of our problem is a non-injective version of the answer given by Ferrand \cite{Ferrand} to a conjecture of Lichnerowicz. 

\begin{theorem}[Lichnerowicz conjecture]\label{lich}\index{Lichnerowicz!conjecture}  Let $K<\infty$.
Up to quasiconformal equivalence, the only compact $n$-manifold which admits a non-compact family of $K$-quasiconformal mappings is the standard $n$-sphere ${\Bbb S}^n$
\end{theorem}

A noninjective rational map $f $ of a closed manifold $M$ will have the semi-group $\{f^n\}_{n=1}^{\infty}$  non-compact - the Julia set of $f$ is always non-empty. Lelong's result suggests the existence of such a map should imply severe restrictions on the manifold $M$. The first of these is the following obstruction for the existence, \cite{MMP}.
\begin{theorem} If $M^n$ is a closed $n$-manifold and $f : M^n \to M^n$ is a non-injective quasiregular rational mapping, then there exists a non-constant quasiregular mapping $g : {\Bbb R}^n \to M^n$.
\end{theorem}
This theorem is proved using a version of Zalcman's lemma.  The semigroup $\{f^{\circ n}\}$ generated by $f$ cannot be normal as the degree increases.  Therefore there is a point $x_0\in M^n$ where the iterates fail to be normal.  Choose a quasiconformal chart $\varphi:{\Bbb B}^n\to M^n$ with $\varphi(0)=x_0$.  The proof is then in showing one can balance the dynamics in $M^n$ and the scaling in ${\Bbb R}^n$ through a judicious choice of sequence $\lambda_n\to 0$ so that
\[  f^{\circ n}(\varphi(\lambda_n x)) \]
is normal and converges to $g$ as $n\to\infty$.  The limit is obviously defined on ${\Bbb R}^n$.

\medskip

Manifolds admitting such a map $g$ are called {\em qr-elliptic},\index{qr-elliptic} and answering a question of Gromov, Varopoulos, Saloff-Coste and Coulhon \cite{VSCC}  showed that such manifolds must in turn have a fundamental group of at most polynomial growth.

\begin{corollary} If $M^n$ is a closed $n$-manifold and $f : M^n \to M^n$ is a non-injective rational quasiregular mapping, then $\pi_1(M^n)$ has polynomial growth.  Hence $M^n$ cannot admit a metric of negative curvature.
\end{corollary}

In fact there are further important consequences of qr-ellipticity.  Bonk and  Heinonen established an upper bound on the dimension of the de Rham coholomogy ring $H^*(M)$ for any closed oriented qr-elliptic $n$-manifold \cite{BH2} and this was generalised by Pankka in other directions \cite{Pankka} to consider mappings whose distortion function is bounded in $L^p$, \cite{Pankka}.  This has obvious implications in generalising this corollary.  Rickman has shown  that  the $4$-manifold $M=({\Bbb S}^2\times {\Bbb S}^2)\# ({\Bbb S}^2\times {\Bbb S}^2)$ is qr-elliptic,  this gives an example which is  nontrivial, simply connected and closed  \cite{rickqr}.  We do not know if this manifold admits a rational quasiregular mapping though. 

\medskip

The generalized Lichnerowicz problem\index{Lichnerowicz!problem} seeks to determine  all closed manifolds which admit non-injective rational quasiregular mappings:  From  \cite{MMP} we have the following.

\begin{theorem} Let f be a non-injective rational quasiregular map of the closed manifold $M$ and suppose that $f$ is locally homeomorphic, so the branch set $B_f = \emptyset$. Then $M $ is the quasiconformal image of a Euclidean space form. \end{theorem}

\medskip

\noindent By a Euclidean space form we mean the quotient of ${\Bbb R}^n$ under a Bieberbach group (co-finite volume lattice) ${\Gamma} \subset Isom^+({\Bbb R}^n)$. The two other types of space forms are the quotients by torsion free co-finite volume lattices of isometries of the $n$-sphere and of hyperbolic $n$-space.  As a sort of converse we also have the following.
\begin{theorem} If $M$ is quasiconformally equivalent to a Euclidean space form, then M admits no branched quasiregular (and in particular no branched quasiregular rational) mappings. \end{theorem}
 
In the case of the sphere, lens spaces and other spherical manifolds the existence of rational quasiregular maps is due to \cite{IM1} and  Peltonen \cite{Pel}. These results suggest that there are few such mappings in three or more dimensions as compared with the space of rational functions of the Riemann sphere ${\Bbb S}^2={\hat{{\Bbb C}}}$.

\begin{theorem} Any non-injective rational quasiregular map of a closed Euclidean space form $M$ is the quasiconformal conjugate of a conformal map. \end{theorem}

We remark that, in this second result, we no longer suppose that the map has to be locally injective. This result is   surprising  because it is false for globally injective mappings. Indeed,  Mayer shows that there are uniformly quasiconformal (even bi-Lipschitz) maps of three (or higher) dimensional tori which cannot be quasiconformally conjugate to a conformal map \cite{Mayer2}.
Next, we can distinguish space forms according to the type of rational quasiregular maps they support:

\begin{theorem}   If $M$ is a closed space form, then we have the following characterization:
\begin{enumerate}
\item $M$ admits a branched quasiregular rational map if and only if $M$ is a spherical space form. 
\item $M$ admits a non-injective, locally injective quasiregular rational map if and only if $M$ is a Euclidean space form.
\item $ M$ admits no non-injective quasiregular rational map if and only if $M$ is a hyperbolic space
form.
\end{enumerate}
\end{theorem}

\section{Quasiconformal Group Actions} 

This final section is intended as an introduction to the theory of quasiconformal groups.
The theory is modelled on the classical theory of discrete groups of M\"obius
transformations and the connections with
and hyperbolic geometry and low-dimensional topology are well known.  

\medskip

A group ${\Gamma} $ of self-homeomorphism of a
domain $\Omega\subset \hat {{\Bbb R}}^n$ is called a {\em quasiconformal group}\index{quasiconformal! group} if there some $K <\infty$ such that each element of ${\Gamma} $ is $K$--quasiconformal.   

\medskip

As a consequence of Theorem \ref{qcgroup} we know that every quasiconformal group\index{quasiconformal!group} admits an invariant conformal structure.  Thus quasiconformal groups are really the groups of conformal transformations of (bounded measurable) conformal structures.

\medskip

Quasiconformal groups
were introduced by Gehring and Palka \cite{GP}  in their study of quasiconformally
homogeneous domains. They asked whether in fact every quasiconformal group is the
quasiconformal conjugate of a M\"obius group.  Sullivan and Tukia established this in two
dimensions,  see \cite{Tukia2D}.   In
higher dimensions the first example of a quasiconformal group not conjugate to a
M\"obius group was given by Tukia
\cite{TukiaNM}. He gave an example of a connected Lie group acting quasiconformally on
${\Bbb R}^n$,  $n\geq 3$,  which was not isomorphic to a M\"obius group.  Tukia's group was
in fact constructed as the topological conjugate of a Lie group.  The obstruction to quasiconformal
conjugacy to a M\"obius group lies in the fact that  the orbit of a point under the group was constructed to be the product
of an infinite Von-Kock snowflake and
${\Bbb R}^{n-2}$.  This orbit is certainly not quasiconformally equivalent to a hyperplane,  the orbit of a point
in an isomorphic M\"obius group.  Certain discrete subgroups
of Tukia's group were also shown not to be quasiconformally conjugate to M\"obius groups
\cite{MartinNM}.  Generalising these examples,  McKemie 
\cite{McKemie} has shown that similar examples both in the discrete and non-discrete case
can be found with $K$ arbitrarily close to $1$.  

\medskip

Examples from topology of ``exotic'' smooth involutions also give finite quasiconformal
groups not conjugate to M\"obius groups.   For example Giffen  \cite{Giffen} shows how
to construct a smooth periodic transformations of the $n$--sphere,  $n\geq 4$,  whose
fixed point set is a knotted co-dimension $2$--sphere.  The fixed point set of any
M\"obius transformation, or its topological conjugate,  must be unknotted. Clearly any
finite group of diffeomorphisms is a quasiconformal group.   

A further important example of an infinite quasiconformal group not topologically
conjugate to a M\"obius group was given in $3$--dimensions by Freedman and Skora \cite{FS}.  This example differs from the others in that the topological fact used 
concerns the linking of the fixed point sets of elliptic elements.  Such linking can be
more complicated for quasiconformal groups than for M\"obius groups.

There are important applications of quasiconformal groups.  Even in $1$--dimension (where
the term {\em quasisymmetric} group is used).  Here Gabai \cite{Gabai} and Casson--Jungreis \cite{CJ}
proved,  building on important earlier work of Tukia \cite{Tukia1D}, Zieschang \index{Zieschang} and others,  that discrete
quasisymmetric groups are topologically conjugate to M\"obius groups.  This is a far reaching generalisation of the Nielsen realisation problem. It was already known from work of Mess and Scott,  that this would also imply an important result in low-dimensional
topology,  namely the Seifert fibered space conjecture: a compact $3$--manifold with
a normal infinite cyclic subgroup of its fundamental group is a Seifert fibered space.
 
\medskip 
In this brief survey we will only outline the basic facts such as the  classical trichotomy classification of elements into  elliptic, parabolic
and loxodromic.   We then recall two fundamental results in the area.  The first asserts that a
``sufficiently large'' discrete quasiconformal group is the quasiconformal conjugate of a
M\"obius group.  Secondly, we show that quasiconformal groups are Lie
groups;  the quasiconformal version of the  Hilbert--Smith conjecture.

\subsection{Convergence Properties}

Suppose that ${\Gamma} $ is a quasiconformal group of self-homeomorphisms of a domain
$\Omega\subset
{\Bbb R}^n$.  Then  ${\Gamma} $ is {\em discrete}\index{discrete quasiconformal group} if it contains no infinite sequence of elements
converging locally uniformly in
$\Omega$ to the identity.  Since the identity is isolated in a discrete group ${\Gamma} $ it follows that 
each element of
${\Gamma} $ is isolated in
${\Gamma} $ in the compact open topology.  The group
${\Gamma} $ is said to be {\em discontinuous}\index{discontinuous} at a point $x\in \Omega$ if there is a neighbourhood
$U$ of $x$ such that for all but finitely many
$g\in {\Gamma} $ we have $g(U) \cap U=\emptyset$.  We denote by $O({\Gamma} )$ the set of all $x\in
\Omega$ such that ${\Gamma} $ is discontinuous at $x$.  $O({\Gamma} )$ is called the {\em ordinary set}\index{ordinary set} of
${\Gamma} $.  The set
\begin{equation}
L({\Gamma} ) = \Omega \setminus O({\Gamma} )
\end{equation}
is called the {\em limit set}\index{limit set} of ${\Gamma} $.  Clearly $O({\Gamma} )$ is open and $L({\Gamma} )$ is closed.  Both
sets are ${\Gamma} $--{\em invariant}.  That is 
\begin{equation}
g(O({\Gamma} ))=O({\Gamma} ) \hskip30pt {\rm and} \hskip30pt g(L({\Gamma} ))=L({\Gamma} ),
\end{equation}
for each $g\in {\Gamma} $.
A discontinuous group (one for which $O({\Gamma} )\neq \emptyset$) is discrete,  the converse is false. 

 \medskip

The following theorem,  a consequence of the compactness results we discussed earlier,  is central to what follows.  
 
\begin{theorem}\label{lemma:discrete} If ${\Gamma} $ is a
discrete quasiconformal group and $\{g_j\}$ is an infinite sequence in ${\Gamma} $,  then there are
points
$x_0$ and $y_0$ and a subsequence $\{g_{j_k}\}\subset\{g_{j}\}$ for which we have $g_{j_k}\to x_0$ locally uniformly in
$\hat {{\Bbb R}}^n\setminus\{y_0\}$ and 
$g_{j_k}\to y_0$ locally uniformly in $\hat {{\Bbb R}}^n\setminus\{x_0\}$
\end{theorem}
It is the above compactness property,   dubbed the {\em convergence property},\index{convergence property} which led
to the theory of convergence groups, see \cite{GM,M3,Tukia1D} which are closely related to
Gromov's  theory of hyperbolic groups.\index{convergence group}

\subsection{The Elementary Quasiconformal Groups}

A discrete quasiconformal group ${\Gamma} $ is said to be {\em elementary}\index{elementary} if $L({\Gamma} )$ consists of
fewer than three points.  For each $g\in {\Gamma} $,  a quasiconformal group,  we set
\begin{eqnarray*}
{\rm ord}(g) & = & \inf \{m>0 : g^m=identity\} \\ {\rm fix}(g) & = &
\{x\in\hat {{\Bbb R}}^n : g(x)=x \}.
\end{eqnarray*}
Let ${\Gamma} $ be a discrete quasiconformal group.  It is an easy consequence of the
convergence properties that if $g\in {\Gamma} $ and ${\rm ord}(g)=\infty$,  then $1 \leq \#{\rm fix}(g) \leq 2$ . 

\medskip

We define three types of elements in a discrete quasiconformal group as
\begin{itemize}
\item $g$ is {\em elliptic};\index{quasiconformal!elliptic} if ${\rm ord}(g)<\infty$;
\item $g$ is {\em parabolic};\index{quasiconformal!parabolic} if ${\rm ord}(g)=\infty$ and $\#{\rm fix}(g)=1$;
\item $g$ is {\em loxodromic};\index{quasiconformal!loxodromic} if ${\rm ord}(g)=\infty$ and $\#{\rm fix}(g)=2$.
\end{itemize}
In a discrete quasiconformal group   this list of elements is exhaustive. 
However, in the non-discrete case there is one other type of element which  needs
to be considered.  If  ${\Gamma} $ is a quasiconformal group and $g\in {\Gamma} $ has the
property that there is a sequence of integers $k_j\to\infty$ for which
\[ g^{k_j}\to {\rm identity} \]
uniformly in $\hat {{\Bbb R}}^n$,   then we call $g$ a {\em quasirotation}.  In this case $\langle g \rangle$ has a nice structure,  we shall see in a moment that its closure in the space of homeomorphisms is a compact abelian Lie group.
 \begin{theorem}  Let  ${\Gamma} $ be a discrete quasiconformal group.  
\begin{itemize} 
\item The limit set  $L({\Gamma} )=\emptyset$ if and only if ${\Gamma} $ is a finite group of elliptic elements.
\item The limit set consists of
one point, $L({\Gamma} )=\{x_0\}$, if and only if ${\Gamma} $ is an infinite group of elliptic or
parabolic elements.
\item   The limit set consists of
two points, $L({\Gamma} )=\{x_0,y_0\}$, if and only if ${\Gamma} $ is an infinite group of loxodromic
elements which fix $x_0$ and $y_0$ and elliptic elements which either fix or interchange
these points.  In addition ${\Gamma} $ must contain at least one loxodromic element and at most
finitely many elliptic elements which fix $x_0$ and $y_0$.  
\end{itemize}
\end{theorem} 

This leads to the classification.

\begin{theorem}\label{theorem:classification}  If ${\Gamma} $ is a discrete quasiconformal group, 
then each $g\in {\Gamma} $ is either elliptic, parabolic or loxodromic.  Moreover,  $g$ and $g^k$
are always elements of the same type for each integer $k\neq 0$.
\end{theorem}
 \subsection{Non-elementary quasiconformal groups and the conjugacy problem.}

In this section we record a few observations about non-elementary groups.  
We first note the following.
\begin{itemize}
\item  If $g$ is a parabolic element with fixed point
$x_0$ in a quasiconformal group
${\Gamma} $, then
\begin{equation}
\lim_{j\to \infty} g^j = x_0 \hskip30pt {\rm and} \hskip30pt \lim_{j\to \infty} g^{-j} =
x_0 
\end{equation}
locally uniformly in $\hat {{\Bbb R}}^n\setminus\{x_0\}$
\item  If $g$ be a loxodromic element with fixed points
$x_0$,
$y_0$ of a quasiconformal group ${\Gamma} $.  Then these points can be labelled so that
\begin{equation}\label{eqn:attr}
\lim_{j\to \infty} g^j = x_0 \hskip30pt {\rm and} \hskip30pt \lim_{j\to \infty} g^{-j} =
y_0 
\end{equation}
locally uniformly in $\hat {{\Bbb R}}^n\setminus\{y_0\}$ and  $\hat {{\Bbb R}}^n\setminus\{x_0\}$ respectively
\end{itemize}  
\subsubsection{The Triple Space}\index{triple space}

A M\"obius group of $\hat {{\Bbb R}}^n$ has the very useful property that it extends to
the upper-half space
\begin{equation}
{\Bbb H}^{n+1}=\{x\in{\Bbb R}^{n+1}:x=(x_1,x_2,\ldots,x_{n+1}),  x_{n+1}>0\}
\end{equation}
via the Poincar\'e extension as a M\"obius group. The existence of such an extension is
unknown for quasiconformal groups and would have important topological consequences.  Next we introduce an alternative for the
upper-half space for which any group of homeomorphisms of $\hat {{\Bbb R}}^n$ naturally extends to. 
This substitute is the {\em triple space}\index{triple space} $T^n$,  a $3n$--manifold defined by
\begin{equation}
T^n =\{(u,v,w):u,v,w\in\hat {{\Bbb R}}^n \mbox{ and $u,v,w$ are distinct }\}.
\end{equation}
There is a natural projection $p:T^n\to {\Bbb H}^{n+1}$ defined by the property that $p(u,v,w)$ is the orthogonal
projection of $w$ (in hyperbolic geometry) onto the hyperbolic line joining $u$ and $v$. 
This map has the property that if $X\subset {\Bbb H}^{n+1}$ is compact, then
$p^{-1}(X)\subset T^n$ is compact.  

Given a self-homeomorphism $f$ of $\hat {{\Bbb R}}^n$,  there is a natural action of $f$ on $T^n$, 
which for notational simplicity we continue to call $f$,  by the rule
\begin{equation}
f(u,v,w) = (f(u),f(v),f(w)), \hskip30pt (u,v,w)\in T^n.
\end{equation} 
If $f$ is in fact a M\"obius transformation of ${\Bbb H}^{n+1}$ we find that the projection $p$
commutes with the action of $f$ on $T^n$,  that is $f\circ p = p\circ f:T^n \to{\Bbb H}^{n+1}$.

Using the convergence properties of quasiconformal
groups it is easy to see  that a quasiconformal group of $\hat {{\Bbb R}}^n$ is discrete if
and only if it acts discontinuously on $T^n$.

\subsubsection{Conjugacy results}

In order to establish the best known result on the quasiconformal conjugacy of a quasiconformal group to a M\"obius groups we need to discuss a special type of limit point.  The following definition is easiest to use in the quasiconformal setting,  though it has useful purely
topological counterparts.

\medskip 
A point
$x_0\in L({\Gamma} )\setminus \{\infty\}$ is called a {\em conical}\index{conical limit point} limit point if there is a sequence of numbers
$\{\alpha_j\}$, 
$\alpha_j \to
0$, and a sequence $\{g_j\}\subset {\Gamma} $ such that the sequence 
\[  h_j(x) = g_j(\alpha_j x + x_0) \]
converges  locally uniformly in ${\Bbb R}^n$ to a quasiconformal mapping $h:{\Bbb R}^n \to \hat {{\Bbb R}}^n\setminus \{y_0\}$.  
We extend the definition to include $\infty$ in
the usual manner. The term {\em radial limit point} is also common in the literature. We should also compare the definition of conical limit point with the
conclusion of the Zalcman Lemma \ref{zalcman}.  The next result is elementary.

\begin{lemma}  Let $x_0$ be a loxodromic fixed point of a $K$--quasiconformal group ${\Gamma}$.  Then $x_0$ is a conical limit point.
\end{lemma}

The loxodromic elements of a quasiconformal group are always quasiconformally
conjugate to M\"obius transformations \cite{GM}.  However, this fact relies on some
quite deep topology.  Parabolic elements are
not always conjugate to M\"obius transformations \cite{Mayer3} and as observed above there are examples of elliptic elements which are not conjugate to
M\"obius transformations.  We now  relate  the regularity of the conformal structure at a conical limit point and the question of conjugacy.  The proof is basically through a linearisation procedure which enables one to move to an invariant conformal structure which is constant after a quasiconformal conjugacy.  This quickly gives us the quasiconformal conjugacy to a M\"obius group.

\begin{theorem}\label{theorem:conjugacy}  Let ${\Gamma} $ be a discrete $K$--quasiconformal group
and let
$G_{\Gamma} $ be a
${\Gamma} $--invariant conformal structure.  Suppose that $G_{\Gamma} $ is continuous in measure at a
conical point of ${\Gamma} $.  Then there is a quasiconformal homeomorphism $h$ of $\hat {{\Bbb R}}^n$ such
that
$h \circ {\Gamma}  \circ h^{-1}$ is a M\"obius group.
\end{theorem}
Since  a measurable map is continuous in measure almost everywhere,  Theorem
\ref{theorem:conjugacy} implies that a quasiconformal group whose conical limit set has
positive measure is the quasiconformal conjugate of a M\"obius group.  There are a number
of results asserting that the set of conical limit points is large,  here are two.

\begin{theorem}\label{conicalqc}  Let ${\Gamma} $ de a discrete quasiconformal group.  Suppose
that either
\begin{itemize}
\item  the action of ${\Gamma} $ on the triple space is cocompact,  that is $T^n/{\Gamma} $ is compact, 
or
\item  the group ${\Gamma} $ can be extended to a quasiconformal group $\tilde{{\Gamma} }$ of ${\Bbb H}^{n+1}$ in such a way
that 
${\Bbb H}^{n+1}/\tilde{{\Gamma} }$ is compact.
\end{itemize}
Then $L({\Gamma} )=\hat {{\Bbb R}}^n$ and every limit point is a conical point.  Thus ${\Gamma} $ is the
quasiconformal conjugate of a M\"obius group.
\end{theorem}
  
 \subsection{Hilbert--Smith Conjecture}

 Hilbert's fifth problem  was formulated following Lie's development of the theory of continuous groups.  It  has been interpreted to ask if every finite-dimensional locally Euclidean topological group is necessarily a Lie
group.  This problem was solved by von Neumann in 1933 for compact groups and by Gleason and Montgomery and
Zippin in 1952 for locally compact groups,  see
\cite{MZ} and the references therein.

A more general version of the fifth problem asserts that among all locally compact groups ${\Gamma}$  only Lie
groups can act effectively on finite-dimensional manifolds.  This problem has come to be called the
Hilbert--Smith Conjecture.  It follows from the work of Newman and of Smith  together with the
structure theory of infinite abelian groups that the conjecture reduces to the special case
when the group ${\Gamma}$ is isomorphic to the $p$-adic integers.  
In 1943 Bochner and Montgomery solved this problem for actions by diffeomorphisms. The Lipschitz
case was establised by Repov\v{s} and \v{S}\v{c}epin \cite{RS}.  In the quasiconformal case we have the following result \cite{Mhil2}.

\begin{theorem}\label{theorem:Hilbert}  Let ${\Gamma}$ be a locally compact group acting 
effectively by quasiconformal homeomorphisms on a Riemannian manifold.  Then ${\Gamma}$ is a Lie
group.
\end{theorem}

Here we wish to make the point that there is no {\em a priori} distortion bounds assumed for elements of ${\Gamma}$.  If one assumes {\em a priori} bounds on the distortion of elements ${\Gamma}$, then
precompactness of the family of all $K$--quasiconformal mappings enables the local
compactness hypothesis in Theorem \ref{theorem:Hilbert} to be dropped.  The hypothesis
of effectiveness (that is the hypothesis that the representation of ${\Gamma}$ in the
appropriate homeomorphism group is faithful) is redundant if we give
${\Gamma}$ the topology it inherits from the compact open topology of maps.  We usually view ${\Gamma}$ simply as a
quasiconformal transformation group.

\begin{corollary}
Let ${\Gamma}$ be a quasiconformal group acting on a Riemannian manifold.  Then
${\Gamma}$ is a Lie group.
\end{corollary}

This result has an important consequence in the property of analytic continuation and also uniqueness statements for solutions of quite general Beltrami systems,  see \cite{Mart1}.

\end{document}